\documentclass[11pt,letterpaper]{article}
\usepackage[dvips,linkcolor=blue,urlcolor=blue,citecolor=blue,linktocpage,bookmarks,pdfstartview=FitBH, breaklinks = true]{hyperref}
\usepackage{fancyhdr}
\usepackage[all]{hypcap}
\usepackage{latexsym}
\usepackage{verbatim}
\usepackage{setspace}
\usepackage[centertags]{amsmath}
\usepackage{amsfonts}
\usepackage{amssymb}
\usepackage{amsthm}
\usepackage{mathrsfs}
\usepackage{newlfont}
\usepackage{natbib}
\usepackage{graphicx}
\usepackage{lscape}
\usepackage{slashbox}
\usepackage{bm}
\usepackage{appendix}
\usepackage{rotating}
\usepackage{url}
\usepackage{float}
\usepackage{enumerate}
\usepackage{datetime}
\usepackage{fancyvrb} 
\usepackage{authblk} 

\numberwithin{equation}{section}

\theoremstyle{plain}
\newtheorem{theorem}{Theorem}[section]

\theoremstyle{definition}

\theoremstyle{remark}

\numberwithin{equation}{section}

\def\0{\mathbf{\dobold 0}}

\def\1{\mathbf{\dobold 1}}

\def\hbbetaS+{\bm{\hat\beta}^{\textrm{S+}}}

\def\D2{\Delta^2} 

\title{Penalty, Shrinkage, and Preliminary Test Estimators under Full Model Hypothesis}

\author[1]{A. K. Md. Ehsanes Saleh}
\author[2]{Enayetur Raheem}
\affil[1]{Carleton University, Ottawa, ON, Canada}
\affil[2]{University of Northern Colorado, Greeley, CO, USA}

\date{}

\begin{document}

\maketitle


\vspace {.5in}

\centerline{\bf Abstract}

\bigskip
In the development of efficient predictive model, the key is to identify the suitable predictors for a given formulated model--may be a multiple linear model or generalized linear models. In search of an efficient technique, several popular estimators appeared in the literature namely, subset selection, non-zero garrott, ridge regression, Bridge estimator, LASSO, adaptive LASSO, SCAD, elastic net among others. This paper considers a multiple regression model and compares, under full model hypothesis, analytically as well as by simulation methodology, the performance characteristics of some popular penalty estimators such as ``ridge regression'' (RR), LASSO, adaptive LASSO, SCAD, and ``elastic net'' (EN) versus Least Squares Estimator (LSE), restricted estimator (RE), preliminary test estimator (PTE), and Stein-type estimators (SE and PRSE) when the dimension ($p$) of the parameter space is smaller than the sample space dimension ($n$). We find that RR uniformly dominates LSE, RE, PTE, SE and PRSE while LASSO, aLASSO, SCAD, and EN uniformly dominates LSE only. Further, it is observed that neither penalty estimators nor Stein-type estimator dominate one another. However, LASSO, aLASSO, SCAD, and EN dominate the Stein-type estimators over a significant portion of the parameter space while LASSO-group uniformly dominates the Stein-type estimators when the dimension of parameter space ($p$) is large. Relative efficiency of EN is a decreasing function of proportion of $L_1$ penalty $(\alpha)$ and correlation $(r)$ among covariates depending on the value of $(\alpha, r)$, and increasing for $\Delta^2>\Delta^2_{(\alpha,r)}$. EN75 dominates both types of Stein-type estimators.  Finally, we observed that RE and RR outperform all estimators at the origin of the parameter space. Our conclusions are based on the analysis of the mean-squared-errors and relative efficiencies with related tables and graphs.



\section{Introduction}\label{sec:intro}

The estimation of parameters of a model with ``uncertain prior information'' on parameters of interest began with \cite{bancroft:1944} in the classical front. But a breakthrough came when \cite{stein:1956} and \cite{james:stein:1961} proved that the sample mean in a multivariate normal model is not admissible under a quadratic loss, for dimension more than two. This very result gave birth to a class of shrinkage estimators of various form and setups. A partial document on preliminary test and Stein-type estimators are given by \cite{judge:bock:1978}. The Stein-type estimators have been reformulated and expanded by \citet[Ch 4.4.3]{saleh:2006} which includes asymptotic and nonparametric methods. Due to the immense impact of Stein's results, scores of technical papers appeared in the literature covering various areas of applications. In 1970, \citeauthor{hoerl:kennard:1970} introduced the ``ridge regression'' estimators which opened the door for a new class of shrinkage estimators known as ``penalty estimators'' based on ``Tikhonov regularization.'' Ridge regression (RR) methodology is a minimization of the least-squares criterion subject to $L_2$ penalty. Ridge regression combats multicollinearity problem in linear models, and is the precursor of the problem of estimation and variable selection. Thus, the estimation subject to penalty function was born. By minimizing least-square criterion subject to $L_p$ penalty function, \cite{frank:friedman:1993} defined a class of ``Bridge estimators'' generalizing the ridge regression estimators. The least-squares criterion with $L_p$-penalty may be written as 
\[
(Y-X\beta)'(Y-X\beta) + \lambda_n \bm{1}'_p|\beta|^\gamma, |\beta|^\gamma =(|\beta_1|^\gamma, \ldots, |\beta_p|^\gamma)
\]
For $\gamma=2$, we get ridge estimates, and $\gamma=1$ relates to LASSO, introduced by \cite{tibshirani:1996}, which has become very popular and an exciting penalty estimator. This estimator is related to the estimators such as ``non-negative'' garotte by \cite{breiman:1996}, smoothly clipped absolute deviation (SCAD) by \cite{fan:li:2001}, elastic net by \cite{zou:hastie:2005}, adaptive LASSO (aLASSO) by \cite{zou:2006}, hard threshold LASSO by
\cite{Belloni:Chernozhukov:2013}, and many other versions. 

This paper is devoted to the study of the performance characteristics of the penalty estimators like ``ridge'', LASSO, aLASSO, SCAD, and elastic net versus Stein-type and preliminary test estimators (PTE). Instead of sub-hypothesis approach, we consider full hypothesis approach to consider the comparisons of penalty and Stein-type estimators. 
In this respect, the paper by \cite{Draper:Craig:1979} is informative with respect to ``ridge'' and Stein-type estimators. An important characteristic of LASSO-type estimators is that they provide simultaneous estimation and selection of coefficients in  linear models and can be applied when dimension of the parameter space exceeds the dimension of the sample space. This paper points to the useful aspects of LASSO-group and ridge estimators as well as the limitations thereof. Conclusions are obtained based on simulated mean-squared errors (MSE) and relative efficiency tables and graphs. 

The organization of the paper is as follows. Section 2 discuss the penalty, shrinkage and preliminary test estimators. Section 3 contains details of analysis of relative efficiencies of estimators with tables and graphs of findings. Conclusions are provided in Section 4.

\section{Linear Model and Estimators}\label{sec:two}
Consider the linear multiple regression model 
\begin{equation}\label{eq2.1}
Y=X\beta + e
\end{equation}
where $Y=(y_1, y_2, \ldots, y_n)',$ $\beta=(\beta_1, \beta_2, \ldots, \beta_p)',$ $X$ is the $n \times p$  design matrix and 
$e=(e_1, e_2, \ldots, e_n)'$ iid error $n$-vector. Further we assume that $E(e)=0$ and $E(e'e)=\sigma^2 I_n$.

It is well-known that the LSE of $\beta$ is given by 
\begin{equation}\label{eq2.3}
\tilde{\beta}_n = (X'X)^{-1}X'Y
\end{equation}
which we designate as the ``unrestricted'' LSE of $\beta$. LSE $\tilde{\beta}_n$ is the ``best linear unbiased estimator (BLUE)'' of $\beta$. Next, we consider few penalty estimators.

Penalty estimators belong to a class of restricted estimators. Simplest example is that of estimators belonging to a linear subspace, namely, $\beta$ coefficients belonging to the subspace defined by $H\beta=h$, where $H$ is a $q \times p$ matrix and $h$ is a $q$-vector of real numbers; $\beta$ being the $p$-vector of coefficients. 

The estimators of $\beta$ is obtained by minimizing the LS criterion subject to $H\beta=h$. Explicitly we may write 
\[
\mathop{\min}_{\beta \in R^p} (Y-X\beta)'(Y-X\beta) + \lambda(H\beta-h)
\]
The solution for this problem is the ``restricted estimator'' and the tuning parameter can be explicitly obtained giving
\begin{equation}\label{eq2.4}
\hat{\beta}_n = \tilde{\beta}_n -C^{-1}H'(HC^{-1}H')^{-1}(H\tilde{\beta}_n-h), \quad C=X'X
\end{equation}
where $\tilde{\beta}_n = (X'X)^{-1}X'Y$, the LSE.

Instead of linear subspace, one may choose a $p$-sphere defined by $\beta'\beta \leq m^2$ and minimize the LS criterion as
\[
\mathop{\min}_{\beta \in R^p} (Y-X\beta)'(Y-X\beta) + \kappa(\beta'\beta - m^2)
\]
There is a one-one correspondence between $\kappa$ and $m^2.$

This minimization yields the estimator with the tuning parameter, known as the ``ridge parameter'' as
\begin{equation}
\hat{\beta}_n(\kappa)= (I_p + \kappa C^{-1})^{-1}\tilde{\beta}_n,
\end{equation}
known as the ridge regression estimator defined by \cite{hoerl:kennard:1970}. Ridge estimator combats multicollinearity problem in linear models. 

If we consider $p$-regular polygons $\sum_{j=1}^{p}|\beta_j| < t$ as the restriction, we minimize
\begin{equation}\label{eq:2.6}
\mathop{\min}_{\beta \in R^p} (Y-X\beta)'(Y-X\beta) + \lambda\sum_{j=1}^p |\beta_j|
\end{equation}
yielding the solution as LASSO introduced by \cite{tibshirani:1996}, given by
\begin{equation}\label{eq:2.7}
\hat{\beta}^L_n = \left\{
\begin{array}{c}
\sum_{i=1}^{n} \left(x_{i\kappa}\beta_\kappa - x_{i\kappa}y_i\right)+\frac{\lambda}{2}\mbox{sgn}(\beta_\kappa)=0 \\
\mbox{ and } \beta_{\kappa}=0
\end{array} 
\right.
\end{equation}
When $\frac{1}{n}X'X \rightarrow I_p$, the solution may be written as 
\begin{align}\label{eq:2.8}
\hat{\beta}_n^L &= \left(\hat{\beta}_{1n}^L,\ldots, \hat{\beta}_{pn}^L\right)\\
\nonumber
\hat{\beta}_{jn}^{L} &= \mbox{ sgn}(\tilde{\beta}_{jn})\left(|\tilde{\beta}_{jn}| -\frac{\lambda}{2}\right)^{+}, j=1, \ldots, p
\end{align} where $C=I_p$, $|\tilde{\beta}_n|=\left(|\tilde{\beta}_{1n}|, \ldots, |\tilde{\beta}_{pn}|\right)'.$

Actually \cite{frank:friedman:1993} defined the class ``Bridge estimators'' as
\begin{eqnarray}\label{eq:2.9}
\mathop{\min}_{\beta\in R^p}(Y-X\beta)'(Y-X\beta) + \lambda\sum_{j=1}^{p}|\beta_j|^\gamma
\end{eqnarray}
If $\gamma=1, 2$ the solution reduces to the LASSO and ridge estimators, respectively.

LASSO proposed by \cite{tibshirani:1996} simultaneously estimates and makes selection of variables with appropriate interpretation and viral popularity in applications. For computational solution and methodology see \cite{tibshirani:1996} and \cite{efron:etal:2004}. Later \cite{efron:etal:2004} proposed Least Angle Regression (LAR) which is a stepwise regression, and \cite{friedman:hastie:tibshirani:2010} developed an efficient algorithm for the estimation of a GLM with convex penalty. During the course of development of penalty estimators, \cite{fan:li:2001} defined good penalty functions as the one which yield (i) nearly unbiased estimator when true parameter is large to avoid unnecessary modeling bias, (ii) an estimator which is a threshold rule that sets small estimated coefficients to zero to reduce model complexity, and (iii) the resulting estimator to be continuous in the data to avoid instability in the model prediction. 

The three characteristics of penalty function is called the ``oracle properties'' (for definition, see \cite{fan:li:2001} and \cite{zou:2006}). It is well-known that LASSO may not possess oracle properties. As such, \cite{fan:li:2001} defined an estimation called the ``smoothly clipped absolute deviation (SCAD)'' estimator based on the continuous differential penalty function defined by 
\begin{equation}\label{eq:2.10}
P_{\alpha, \lambda}(\beta) = \lambda\left\{I(|\beta| \leq \lambda) + \frac{(\alpha \lambda -|x|)^{+}}{(\alpha-1)} I(|\beta| > \lambda)\right\}, \beta \ge 0
\end{equation}
Here $\lambda>0$ and $\alpha>2$ are tuning parameters. For $\alpha=\infty,$ the expression (\ref{eq:2.10}) is equivalent to $L_1$-penalty. Thus, the SCAD estimator is obtained by 
\begin{equation}\label{eq:2.11}
\hat{\beta}_n^{SCAD} = \mathop{\min}_{\beta \in R^p} \left[(Y-X\beta)'(Y-X\beta) + \lambda \sum_{j=1}^{p}P_{\alpha, \lambda}(\beta_j) |\beta_j|\right]
\end{equation}
Later, \cite{zou:2006} modified the LASSO penalty by using weighted LASSO and defined the ``adaptive LASSO" estimator given by
\begin{equation}\label{eq:2.12}
\hat{\beta}_n^{aL} =  \mbox{arg}\min\left[(Y-X\beta)'(Y-X\beta) + \lambda
_n \hat{W}|\beta|\right]
\end{equation}
where $\hat{W}=(\hat{w}_1, \ldots, \hat{w}_p)'$, $\hat{w}_j=|\beta_{nj}^{*}|^{-\gamma}, \gamma >0$, $|\beta| =(|\beta_1|, \ldots, |\beta_p|)'$ and $\beta_{n}^{*}$ is $\sqrt{n}$-consistent vector for $\beta.$ The equation (\ref{eq:2.12}) is a convex optimization problem and its global minimizer can be solved efficiently (see \cite{zou:2006}). Detailed computational methodology is given by \cite{zou:2006}.

\cite{zou:hastie:2005} proposed ``elastic net''-estimator given by 
\begin{equation}\label{eq:Elastic:Net}
\hat{\beta}_n(\lambda_n)^{EN}=\mbox{arg}\min\left[(Y-X\beta)'(Y-X\beta) + \lambda_n\left\{\alpha \sum_{i=1}^{p}|\beta_j| + (1-\alpha)\sum_{i=1}^{p}|\beta_j|^2\right\}\right]
\end{equation}
where $\lambda_n$ is the tuning parameter and $\alpha\in (0,1)$, which borrows useful properties of LASSO and ridge regression estimators. The first component of the penalty term in (\ref{eq:Elastic:Net}) is the LASSO penalty while the second component is the ridge regression estimates.  Here, $\alpha$ controls the mixing proportion of LASSO and ridge estimators. The value of $\alpha$ can be set judiciously or based on cross-validation. Thus, elastic net estimator is advantageous when there are highly correlated predictors present in the model and a sparse solution is needed. 

\cite{zou:zhang:2009} proposed adaptive elastic net which borrows the useful features of quadratic regularization and adaptive LASSO shrinkage. They showed, under weak regularity conditions, that the adaptive elastic net has oracle property. Adaptive elastic net performs better than other estimators having oracle property when collinearity is present \citep{zou:zhang:2009}.

\cite{yuan:lin:2006} proposed grouped LASSO which can be applied when selection of variables ``in groups'' are essential. This method performs better than traditional stepwise variable selection methods. However, grouped LASSO is often inefficient and inconsistent in variable selection problems. To overcome the estimation efficiency and inconsistency of variable selection in the original grouped LASSO, \cite{Wang:Leng:2008} extended the grouped LASSO and proposed adaptive grouped LASSO as a tool for shrinkage and variable selection in a grouped manner. Adaptive grouped LASSO is as efficient as other oracle-like estimators \citep{Wang:Leng:2008}.

\subsection{PTE and Stein-type Estimators}
For the linear multiple regression model, $Y=X\beta+e$, if we suspect the full hypothesis to be $\beta=0$ (null-vector), then the restricted estimator (RE) $\hat{\beta}_n=0$ and the test for $\beta=0$, vs $\beta \ne 0$ may be based on the statistic
\begin{subequations}
\begin{eqnarray}
\mathcal{L}_n =& \frac{\tilde{\beta}_n'C\tilde{\beta}_n}{s_e^2}, \quad s_n^2 = (n-p)^{-1}(Y-X\tilde{\beta}_n)'(Y-X\tilde{\beta}_n) \\
\nonumber \mbox{ Under the conditions } \\
\nonumber (1)& \frac{1}{n}X'X \rightarrow C \mbox{ as } n\rightarrow \infty \, (C \mbox{ positive definite }), \mbox{ and} \\ 
(2)& \mathop{\max}_{1 \leq j \le n} \bm{x}_j'(\frac{1}{n}X'X)^{-1}\bm{x}_j \rightarrow 0 \mbox{ as } n \rightarrow \infty
\end{eqnarray}
\end{subequations}
where $\bm{x}_j$ is the $j$th row of $X$, and $\mathcal{L}_n \rightarrow\chi^2_p$ -central chi-square variable with $p$ degrees of freedom (df).

Let $\chi^2_p(\alpha)$ be an upper $\alpha$-level critical value from this null distribution; then we may define the PTE of $\beta$ as
\begin{equation}\label{eq:2.14}
\hat{\beta}_n^{PT} = \tilde{\beta}_n -\tilde{\beta}_n \, I(\chi^2 < \chi^2_p(\alpha))
\end{equation}

The PTE is a discrete variable. As a result some optimality properties when we consider assessing its MSE comparison is lost. We may define a continuous version of PTE as the James-Stein-type estimator given by
\begin{equation}\label{eq:2.15}
\hat{\beta}_n^{S}=\tilde{\beta}_n -(p-2)\tilde{\beta}_n \mathcal{L}_n^{-1}
\end{equation}
Note that we have replaced $I(\chi^2 < \chi^2_p(\alpha))$ by $(p-2)\mathcal{L}_n^{-1}$ in the definition of PTE. However, $\hat{\beta}_n^{S}$ has an inherent problem of changing its sign due to the factor $(1-(p-2)\mathcal{L}_n^{-1})$ which may be larger than 1 in absolute value. If that happens, from applied point of view, its interpretation becomes blurred. Thus, we define another estimator, namely, the positive-rule Stein-type estimator (PRSE) as follows:
\begin{equation}\label{eq:2.16}
\hat{\beta}_n^{S+} = \hat{\beta}_n^{S} \, I(\mathcal{L}_n > (p-2))
\end{equation}
Thus, we have defined five estimators, namely, LSE, RE, PTE, JSE, and PRSE here and four penalty estimators namely, Ridge regression (RR), LASSO, aLASSO, and SCAD. We like to compare the performance characteristics of these two groups of estimators for practical utility. 

Next, we define an improved preliminary test (IPT) estimator defined by 
\begin{equation}\label{eq:2.17}
\hat{\beta}_n^{IPT} = \hat{\beta}_n^{PT}\, \left(1-(p-2)\mathcal{L}_n^{-1} \right)
\end{equation}

\subsection{Quadratic Risk-functions for LSE, RE, PTE, Stein-type and Improved PTE}
In this section, we consider the asymptotic distributional bias (ADB) and quadratic risk functions of the above estimators.

Since $\beta=0$ is uncertain, we test the hypothesis, $H_0$, based on the statistic
\begin{equation}\label{eq:2.18}
\mathcal{L}_n= \frac{\tilde{\beta}_n'C\tilde{\beta}_n}{s^2_n}, \quad s^2_n = (n-p)^{-1} (Y-X\tilde{\beta}_n)'(Y-X\tilde{\beta}_n)
\end{equation}
This test is consistent and its power function converges to unity as $n\rightarrow \infty$ for fixed alternatives. Thus, we consider the sequences of local alternatives, $\left\{K_{(n)}\right\}$ defined by
\begin{equation}\label{eq:2.19}
K_{(n)}: \beta_{(n)} = \frac{1}{\sqrt{n}}\, \delta, \quad \delta = (\delta_1, \ldots, \delta_p)'\ne 0
\end{equation}
and the loss function $n(\beta^* - \beta)'\bm{W}(\beta^*-\beta)$. 
For $\delta=0,$ $K_{(n)}\equiv H_0,$ then $\sqrt{n}(\tilde{\beta}_n-\beta_{(n)}) \approx N_p(\bm{0}, \sigma^2\bm{C}^{-1}).$ Hence, for $\bm{W}=I_p$ we have the asymptotic distributional bias, MSE, and risk expressions as
Hence, 
\begin{align}
(1) & ADB(\tilde{\beta}_n) = \bm{0}, \quad ADQR(\tilde{\beta}_n) = \sigma^2\mbox{ tr }C^{-1}\\
(2) & ADB(\hat{\beta}_n)=-\delta, \quad ADQR(\hat{\beta}_n)=\sigma^2 \Delta^2 \\
\nonumber 
(3) & ADB(\hat{\beta}_n^{PT}) = -\delta \mathcal{H}_{p+2} (\chi^2_p(\alpha); \Delta^2) \\
\nonumber 
& ADQR(\hat{\beta}_n^{PT}) = \sigma^2\left[\mbox{tr }C^{-1}\right]\left[1-\mathcal{H}_{p+2}\left(\chi^2_p(\alpha); \Delta^2\right)\right] + \\
& \quad  \sigma^2\Delta^2 \left\{2\mathcal{H}_{p+2}\left(\chi^2_p(\alpha); \Delta^2\right) -\mathcal{H}_{p+4}\left(\chi^2_p(\alpha); \Delta^2\right)\right\}
\end{align}
Similarly,
\begin{align}
\nonumber (4)& ADB(\hat{\beta}_n^{S}) = -(p-2) \delta E\left[\chi^{-2}_{p+2}(\Delta^2)\right] \\
\nonumber
& ADQR(\hat{\beta}_n^{S}) = \sigma^2\left[\mbox{tr }C^{-1}\right]\left[1-(p-2)\left\{2E\left[\chi^{-2}_{p+2}(\Delta^2)\right]-(p-2)E\left[\chi^{-4}_{p+2}(\Delta^2)\right]\right\}\right]\\
& \quad + \sigma^2(p^2-4) \Delta^2 E\left[\chi^{-4}_{p+4}(\Delta^2)\right]
\end{align}

\begin{align}
(5)& ADB(\hat{\beta}_n^{S+}) = ADB(\hat{\beta}_n^{S}) - \delta E\left[\left(1-(p-2)\chi^{-2}_{p+2}(\Delta^2)\right)\, I\left(\chi^2_{p+2}(\Delta^2)<(p-2)\right)\right] \\
\nonumber & ADQR(\hat{\beta}_n^{S+}) = ADQR(\hat{\beta}_n^{S}) \\
\nonumber & \quad - \sigma^2 \left[\mbox{tr }C^{-1}E\left\{\left(1-(p-2)\chi^{-2}_{p+2}(\Delta^2\right)^2 \, I\left(\chi^2_{p+2}(\Delta^2) < p-2\right)\right\}\right] \\
\nonumber 
& \quad + \sigma^2\Delta^2\left\{ 2E\left[\left(1-(p-2)\chi^{-2}_{p+2}(\Delta^2)\right)\, I\left(\chi^{-2}_{p+2}(\Delta^2) < p-2\right)\right] \right.\\
& \quad \left. -E\left[\left(1-(p-2)\chi^{-2}_{p+4} (\Delta^2)\right)^2\, I\left(\chi^2_{p+4}(\Delta^2) < p-2\right)\right] \right\} \\
\nonumber
(6) & ADB(\hat{\beta}_n^{IPT}) = ADB(\hat{\beta}_n^{PT}) \\
& \quad - \delta E\left[\left(1-(p-2)\chi^{-2}_{p+2}(\Delta^2)\right)\,I\left(\chi^2_{p+2}(\Delta^2) <\chi^2_p(\alpha)\right)\right] \\
\nonumber
& ADQR(\hat{\beta}_n^{IPT}) = ADQR(\hat{\beta}_n^{PT}) \\
\nonumber & \quad -\sigma^2 (p-2)\left(\mbox{tr} C^{-1}\right) \left\{ E\left[\chi^{-2}_{p+2}(\Delta^2)\right] + \Delta^2E\left[\chi^{-4}_{p+4}(\Delta^2)\right]\right.\\
\nonumber
 & \quad \left. - E\left[\chi^{-2}_{p+2}(\Delta^2) I\left(\chi^2_{p+2}(\Delta^2)<\chi^2_p(\alpha)\right)\right] \right\}
\\
\nonumber 
& \quad + (p-2) \sigma^2\Delta^2\left\{ (p+2)E\left[\chi^{-4}_{p+4}(\Delta^2)\right] + E\left[\chi^{-2}_{p+4}(\Delta^2) I\left(\chi^2_{p+4}(\Delta^2) <\chi^2_p(\alpha)\right)\right]\right\} \\
(7) & ADB\left(\hat{\beta}_{n(\kappa
)}^{RR}\right) = -\kappa \left(C+\kappa I_p\right)^{-1}\delta\\
& ADQR\left(\hat{\beta}_{n(\kappa)}^{RR}\right) = \sigma^2 \mbox{tr}\left[\left(C+\kappa I_p\right)^{-1} C^{-1}\left(C+\kappa I_p\right)^{-1}\right] + \kappa^2\delta'\left(C+\kappa I_p\right)^{-2}\delta
\end{align}
Here, $\mathcal{H}_{p+2\nu}(\cdot; \Delta^2)$ is the cdf of a noncentral chi-square distribution with $p+2\nu$ d.f. with noncentrality parameter $\Delta^2$, and 

\[
E[\chi^{-2r}_{p+2\nu}(\Delta^2)] = \int_0^{\infty} x^{-2r} d\mathcal{H}_{p+2\nu}(x;\Delta^2)
\]

\subsection{Analysis of Risk Functions of LSE, RE, PTE, SE and PRSE}

\begin{enumerate}[(i)]
\item \emph{Comparison of LSE and RE}: The risk-difference of LSE and RE is given by 
\[
\sigma^2 \left[\mbox{tr } C^{-1} - \Delta^2\right].
\] 
If $\Delta^2<\mbox{tr }C^{-1}$, then RE is better than the LSE, and if $\Delta^2 >\mbox{tr }C^{-1}$ then LSE is better than RE.

\item \emph{Comparison of LSE and PTE}: Here the risk-difference is given by 
\[
\sigma^2\left[\left(\mbox{tr }C^{-1}\right)\mathcal{H}_{p+2} \left(\chi^2_p(\alpha); \Delta^2\right) - \Delta^2 \left\{2\mathcal{H}_{p+2}\left(\chi^2_{p+2}(\alpha); \Delta^2\right)-\mathcal{H}_{p+4}\left(\chi^2_p(\alpha); \Delta^2\right)\right\}\right]
\]
Hence, if 
\[
\Delta^2 \leq \frac{(\mbox{tr }C^{-1})\mathcal{H}_{p+2}\left(\chi^2_p(\alpha);\Delta^2\right)}{\left[2\mathcal{H}_{p+2}\left(\chi^2_p(\alpha); \Delta^2\right)-\mathcal{H}_{p+4}\left(\chi^2_p(\alpha);\Delta^2\right)\right]}.
\]
PTE performs better than the LSE, and LSE performs better than PTE whenever
\[
\Delta^2 > \frac{(\mbox{tr }C^{-1})\mathcal{H}_{p+2}\left(\chi^2_p(\alpha);\Delta^2\right)}{\left[2\mathcal{H}_{p+2}\left(\chi^2_p(\alpha); \Delta^2\right)-\mathcal{H}_{p+4}\left(\chi^2_p(\alpha);\Delta^2\right)\right]}.
\]

\item \emph{Comparison of LSE and Stein-type estimators}: The risk-difference of LSE and James-Stein estimator is given by
\begin{align*}
& \sigma^2 (p-2)\left[\mbox{tr }C^{-1} \left\{2E\left[\chi^{-2}_{p+2}(\Delta^2)\right] - (p+2)E\left[\chi^{-4}_{p+2}(\Delta^2)\right]\right\} -(p+2)\frac{\delta'\delta}{\sigma^2}E\left[\chi^{-4}_{p+4}(\Delta^2)\right]\right]\\
=& \sigma^2(p-2)\mbox{tr }C^{-1}\left[(p-2)E\left[\chi^{-4}_{p+2}(\Delta^2)\right]+\left(1-\frac{(p+2)\delta'\delta}{2\Delta^2\mbox{tr}C^{-1}}\right)2 \Delta^2E\left[\chi^{-4}_{p+4}(\Delta^2)\right]\right]\\
\ge & 0
\end{align*}
whenever 
\[
\frac{\mbox{tr }C^{-1}}{\mbox{Ch}_{\mbox{max}}(C^{-1})} \ge \frac{p+2}{2}
\]
Hence,
\[
ADQR(\hat{\beta}_n^S) \leq ADQR(\tilde{\beta}_n) \quad \forall \Delta^2
\]

Next, we consider the risk-difference of $\hat{\beta}_n^S$ and $\hat{\beta}_n^{S+}$, which is given by 
\begin{align*}
& \sigma^2 \left[\left(\mbox{tr }C^{-1}\right) E\left\{\left(1-(p-2)\chi^{-2}_{p+2}(\Delta^2)\right)^2 I\left(\chi^2_{p+2}(\Delta^2) < p-2\right)\right\}\right]\\
& \quad -\delta'\delta \left\{2E\left[\left(1-(p-2) \chi^{-2}_{p+2}(\Delta^2)\right)I\left(\chi^2_{p+2}(\Delta^2) < p-2\right)\right]\right. \\ 
& \quad - \left. 2E\left[\left(1-(p-2)\chi^{-2}_{p+4}(\Delta^2)\right)^2 I\left(\chi^{-2}_{p+4}(\Delta^2)<p-2\right)\right]\right\}\\
=& \sigma^2\left[\left(\mbox{tr }C^{-1}\right)E\left[\left(1-(p-2)\chi^{-2}_{p+2}(\Delta^2)\right)^2 I\left(\chi^2_{p+2}(\Delta^2) <p-2\right)\right]\right. \\
& \quad \left. + \delta'\delta \left\{E\left[\left((p-2)\chi^{-2}_{p+2}(\Delta^2)-1\right) I\left(\chi^2_{p+2}(\Delta^2)<p-2\right)\right]\right.\right. \\
& \quad \left. + E\left[\left(1-(p-2)\chi^{-2}_{p+4}(\Delta^2)\right)^2 I\left(\chi^2_{p+4}(\Delta^2)<p-2\right)\right]\right\}\\
& \ge 0
\end{align*}
Hence, 
\[
ADQR(\hat{\beta}_n^{S+}) \leq ADQR(\hat{\beta}_n^{S}) \leq ADQR(\tilde{\beta}_n) \quad \forall \Delta^2
\]

The optimum tuning parameter of the ridge estimator is obtained as follows. The AQDR of $\hat{\beta}^{RR}$ is given by 
\[
\sigma^2\mbox{ tr}\left[\left(C+\kappa I_p\right)^{-1}C^{-1}\left(C+\kappa I_p\right)^{-1}\right] + \kappa^2\delta'\left(C+\kappa I_p\right)^{-2}\delta
\]
Now, since $C$ is a positive semidefinite matrix there exists an orthogonal matrix $\Gamma$ such that $\Gamma'C\Gamma=I_p$, then the eigenvalues of $C+\kappa I_p$ are $\{(1+\kappa), \ldots, (1+\kappa)\}$ and we have 
\[
\mbox{tr}\left[\left(C+\kappa I_p\right)^{-1}C^{-1} \left(C+\kappa I_p\right)\right] = \sum \frac{1}{(1+\kappa)^2} = \frac{p}{(1+\kappa)^2}
\]
and
\[
\delta'(C+\kappa I_p)^{-2}\delta = \frac{\alpha'\alpha}{(1+\kappa)^2}, \quad \alpha=\Gamma\delta'.
\]
Hence, ADQR of ridge estimate is given by 
\[
\frac{\sigma^2p}{(1+\kappa)^2} + \frac{\kappa^2 \alpha'\alpha}{(1+\kappa)^2}
\]
And
\begin{align*}
& \frac{\partial}{\partial \kappa} \left\{\frac{\sigma^2p}{(1+\kappa)^2}+\frac{\kappa^2 \alpha'\alpha}{(1+\kappa)^2}\right\} \\
&= 2 \, \frac{1}{(1+\kappa)^3}\, (\kappa\alpha'\alpha - p\sigma^2)
\end{align*}
Therefore, 
\[
\kappa = \frac{p\sigma^2}{\alpha'\alpha} = \frac{p}{\Delta^2}
\]
Hence, $\hat{\beta}_n^{RR}$ dominates $\tilde{\beta}_n$ uniformly for $\kappa \in \left(0, \frac{p}{\Delta^2}\right]$.

\begin{theorem}
$\hat{\beta}_n^{IPT}$ dominates $\hat{\beta}_n^{PT}$ in $(\alpha, \Delta^2) \in (0, 1)\cup (0, \infty)$ for $p\ge 3$.
\end{theorem}
\proof For $p\ge 3$, risk-difference is given by 
\begin{align*}
&\sigma^2(p-2)(\mbox{tr }C^{-1})\left\{E\left[\chi^{-2}_{p+2}(\Delta^2)\right] + \Delta^2E\left[\chi^{-4}_{p+4}(\Delta^2)\right] \right. \\
& \quad \left.-E\left[\chi^{-2}_{p+2}(\Delta^2) I\left(\chi^2_{p+2}(\Delta^2)< \chi^2_p(\alpha)\right)\right] \right\}\\
& \quad -(p-2) \beta'\beta\left\{(p+2)E\left[\chi^{-4}_{p+4}(\Delta^2) I \left(\chi^2_{p+4}(\Delta^2)<\chi^2_p(\alpha)\right)\right]\right\}\\
& \ge 0.
\end{align*}
Hence, $\hat{\beta}_n^{IPT}$ dominates $\hat{\beta}_n^{PT}$ uniformly in ($\alpha, \Delta^2$). If $\chi^2_p(\alpha) < p-2$, then $\hat{\beta}_n^{IPT}$ behaves like $\hat{\beta}_n^{S+}$ and if $\chi^2_p(\alpha) \ge p-2$, then the above dominance result holds.
\end{enumerate}

\section{Simulation} \label{sec:simulation}

We conduct Monte Carlo simulation experiments to study the performance of the Stein-type, pretest, and penalty estimators. In particular, we compare relative efficiencies of the estimators compared to the least squares estimator (LSE). In the simulation studies, mean squared errors (MSE) were computed for each of the estimators, and their relative efficiencies were obtained by taking the ratio of MSE of the estimators to the MSE of LSE. \cite{raheem:ahmed:doksum:2012} have studied through simulation where a sub-hypothesis was tested and relative efficiencies of various shrinkage and penalty estimators were studied in a partially linear regression setup. In this study, we are interested in testing the hypothesis $H_0:\beta=0,$ against the alternative $H_a: \beta \ne 0$. The simulation setup is discussed in the following section. 


We generate $x$-matrix from a multivariate normal distribution with mean vector $\bm{\mu} =\bm{0}$ and covariance matrix $\bm{\Sigma}$. The off-diagonal elements of the covariance matrix are considered to be equal to $r$ with $r=0, 0.2, 0.9$. We consider $n=100$ and various $p$ ranging from 10 to 95 depending on the comparative studies we have performed.

In our setup, $\beta$ is a $p$-vector and a function of $\Delta^2$. When $\Delta^2=0$, $\beta$ is the null vector. $\Delta^2 >0$ is equivalent to ``violation'' of the null hypothesis. We considered 23 different value for $\Delta^2$, which are 0, 0.1, 0.2, 0.3, 0.4, 0.5, 0.6, 0.7, 0.8, 0.9, 1, 1.5, 2, 3, 5, 10, 15, 20, 25, 30, 35, 40, and 50. The way the $\beta$ vector is defined in our setup, a $\Delta^2=0$ indicates that data are generated under null hypothesis, whereas $\Delta^2 >0$ indicates a data set generated under alternative hypothesis. Each realization was repeated 2000 times to obtain bias-squared and variance of the estimated regression parameters. 

Finally, MSEs are calculated for the LSE, RE, Stein-type (S) and positive-rule estimator (S+), Pretest, ridge  regression (RR), and the penalty estimators. The responses were simulated from the following model:
\[
y_i = \sum_{i=1}^{p}x_i\beta_i + e_i
\]
where $e_i \sim N(0, 5^2)$.

We use the following formula to calculate relative efficiency:
\begin{equation}\label{eq:relative:efficiency}
\mbox{Relative Efficieicy} = \frac{\mbox{MSE}(\hat{\beta}_n^{\textrm{LSE}})}{\mbox{MSE}(\hat{\beta}_n^{*})}
\end{equation}
where $\hat{\beta}_n^{*}$ is one of the estimators whose relative efficiency is to be computed.

For comparing the relative efficiencies of the penalty estimators, the data generation setup was slightly modified to accommodate the number of non-zero $\beta$s in the model. In particular, we partitioned $\beta$ as $\beta=(k, q)'$ where $k$ indicates number of nonzero $\beta$s and $q$ indicates $p-k$ zeros--a function of $\Delta^2$. To translate the above, when $p=10,$ and $k=5$, we would have $\beta=(1, 1, 1, 1, 1, 0,0,0,0,0)'$, and the previously mentioned procedure would be used to generate the data. We used R statistical software \citep{r:2014} to carry out the simulation.

In the following, we discuss our results of three different simulation studies to compare (i) James-Stein shrinkage, preliminary tests, ridge, and elastic net (EN) estimators, (ii) penalty estimators (LASSO, aLASSO, and SCAD), and (iii) James-Stein shrinkage and penalty estimators. 

\subsection{Comparison of Shrinkage, Preliminary Tests, Ridge, and Elastic Net Estimators} \label{sec:simulation:compare:JS:Ridge}
In this study, data have been generated with correlation between the $x$'s, $r=0, 0.2, 0.9$ for $n=100$ and $p=10$. In Tables~\ref{tb:rel:eff:cor0} through \ref{tb:rel:eff:cor9}, we present relative efficiencies of restricted estimator (RE), preliminary test estimators (PT) for $\alpha=.05, .10, .15, .25$, improved preliminary test estimator (IPT) for $\alpha=.10$, James-Stein-type shrinkage estimator (S), positive rule shrinkage estimator (S+), the ridge regression (RR) estimator, and elastic net (EN) estimators (EN25, EN50, EN75). Here, EN25 represents elastic net estimate with mixing parameter $\alpha=.25$. Similarly, EN50 and EN75 represent elastic net estimators with $\alpha=.5$ and $\alpha=.75$, respectively. This $\alpha$ is not the same $\alpha$ considered in PT and IPT estimators. Rather, $\alpha$ is the mixing parameter in the elastic net estimator that controls the mixing proportion of ridge and LASSO estimates. We used {\tt glmnet()} R package \citep{glmnet:package} in which the same notation for the mixing parameter was used for elastic net estimation. We retained the notation as is. At the two extremes: $\alpha=0$ gives ridge estimates while $\alpha=1$ gives LASSO estimates. For example, $\alpha=.5$ is a 50-50 mixture of ridge and LASSO. It is observed that for $r=0.9$, EN75 outperforms S and S+ uniformly while EN25 and EN50 outperform S and S+ over a significant portion of the parameter space defined by $\Delta^2$.   See \cite{glmnet:package} and the vignette of the {\tt glmnet()} R-package for more information on elastic net estimation. 

For easier comparison among the estimators we plotted relative efficiencies of RE, PT ($\alpha=.15$ only), S, S+, RR, and EN estimators against $\Delta^2$ in Figures~\ref{fig:ols-ridge-p10-r0}-\ref{fig:ols-ridge-p10-r9}. A horizontal line was drawn at 1 on the $y$-axis to facilitate the comparison among the estimators. For a given estimator, any point above this line indicates superiority of the estimator compared to the LSE in terms of relative efficiency.

The findings of simulation study may be summarized as follows.

\begin{enumerate}[(i)]
\item From Figures~\ref{fig:ols-ridge-p10-r0}-\ref{fig:ols-ridge-p10-r9}, we find that ridge regression outperforms RE, PT, S, and S+ estimators in terms of relative efficiency. The efficiency gain increases as the correlation of the $x$'s increase. 

\item The restricted estimator performs better than the LSE when $\frac{\mbox{tr }C^{-1}}{\Delta^2}\ge 1,$ while LSE dominates when $\frac{\mbox{tr }C^{-1}}{\Delta^2}<1$.  See Figure~\ref{fig:ols-ridge-p10-r9}. 

\item It is well-known that shrinkage and positive-rule shrinkage estimators are always at or better than the LSE. Similarly we find IPT is superior to PTE uniformly in $(\alpha, \Delta^2)$ for $p\ge 3$. Further, when $\Delta^2$ is at or near zero, IPT is the best estimator after ridge in terms of relative efficiency.

\item Neither PTE nor Stein-type estimator dominate one another. 
\end{enumerate}

\subsection{Comparison of Absolute Penalty Estimators}\label{sec:simulation:compare:penalty}

In the second simulation experiment, we compared LASSO, aLASSO and SCAD estimators for $r=0.2, 0.9$ and $p=10, 20, 30$ for varying $\Delta^2$. The results are presented in Tables~\ref{tb:rel:eff:penalty:cor2:p10}-\ref{tb:rel:eff:penalty:p30cor9}. As outlined at the beginning of Section~\ref{sec:simulation}, data were generated for various configurations of $\beta=(k, p-k)'$ where $k$ is the number of 1s and $p-k$ is the number of zeros. 

We find that relative efficiencies tend to decrease as $k$ increases. Next, we observe that relative efficiencies may be ordered as: if $\Delta^2<10$ and $r=0.2$, then $LASSO > aLASSO > SCAD$ and for $\Delta^2 \ge 10$ $aLASSO > LASSO > SCAD$. When $\Delta^2<10$ and $r=0.9$, the order is $SCAD < LASSO < aLASSO$ and for $\Delta^2\ge 10$ $aLASSO > LASSO > SCAD$.

Of the three penalty estimators considered, SCAD was found to be the least efficient in our experiment. We noticed that as $k$ increases, SCAD become less and less efficient compared to LSE. In particular, when $r=.9$, SCAD performs worst compared to LASSO and aLASSO. In other words, SCAD performs relatively better when there are more near-zero coefficients (small $k$) and the correlation among the $x$'s is small. 

\subsection{Comparison of Absolute Penalty and Stein-type Estimators} \label{sec:simulation:compare:penalty:JS}

In this simulation experiment, we compare James-Stein-type shrinkage estimators (S, and S+) with LASSO, aLASSO, and elastic net (EN). We have considered elastic net since it is a mixture of LASSO and ridge estimators. Considering elastic net allows us to eliminate ridge estimator for this part of the simulation studies. This somewhat broadens our scope of comparing various estimators with the limited configurations (i.e, $r$, $p$) that we have considered in our simulation experiments. 

To compare among these estimators, we consider $r=0.2$, $k=0, 1, 3, 5$, $p$ ranging from 10 to 95, and $n=100$. Recall that $k$ is the number $\beta$s equal to 1. Relative efficiencies of LASSO, aLASSO, elastic net, and Stein-type estimators are presented in Table~\ref{tb:penalty:JS:comparision}. Figure~\ref{fig:penalty:JS} graphically shows relative efficiencies of JS-type shrinkage estimator (S), positive-rule shrinkage estimator (PRSE, S+), LASSO (L), aLASSO (aL), and elastic net (EN25, En50, EN75) when compared to the LSE. It is observed for $r=0.9$ that EN75 outperforms S and S+ uniformly while E25 and EN50 outperform S and S+ over a significant portion of the parameter space defined by $\Delta^2$.

Overall, our simulation results show that PRSE (S+) estimator outperforms all other estimators when $k>0$ and when $p\leq 30$. Adaptive LASSO performs best when $k=0$ for all $p$. As $p$ increases LASSO, aLASSO, and elastic net estimators outperform JS-type shrinkage estimators. For $p \ge 60$ and $k>1$, elastic net estimators outperform all other estimators (Table~\ref{tb:penalty:JS:comparision}).


\section{Summary and Conclusion}
\begin{figure}[!h]\label{fig:ols-ridge-p10-r0}
\includegraphics[width=6in]{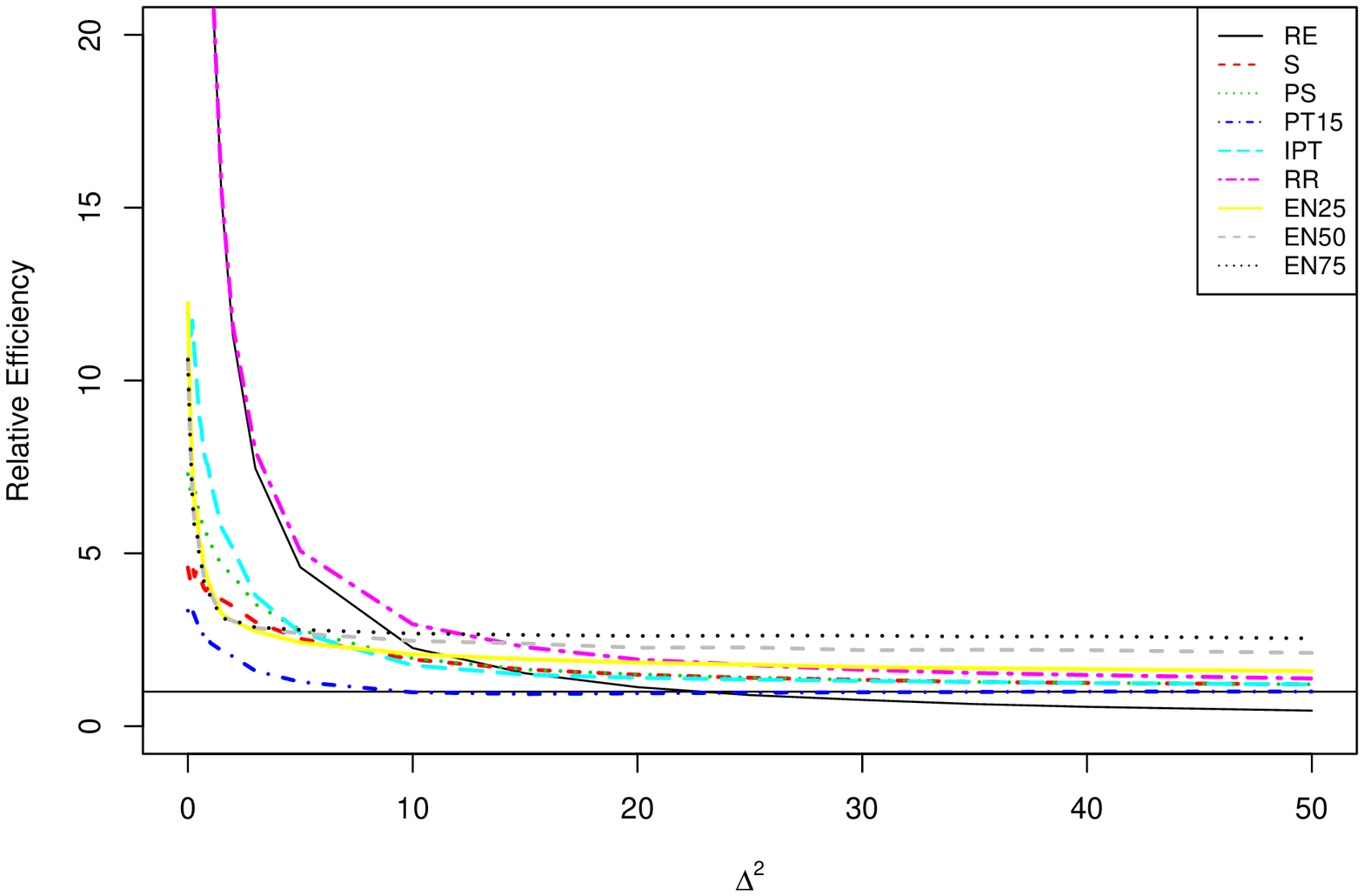}
\caption{Relative Efficiency plots for Shrinkage, Pretest, Ridge, and Elastic Net Estimators when $p=10$, $r=0$, and $n=100$.}
\end{figure}

\begin{figure}[!h]\label{fig:ols-ridge-p10-r2}
\includegraphics[width=6in]{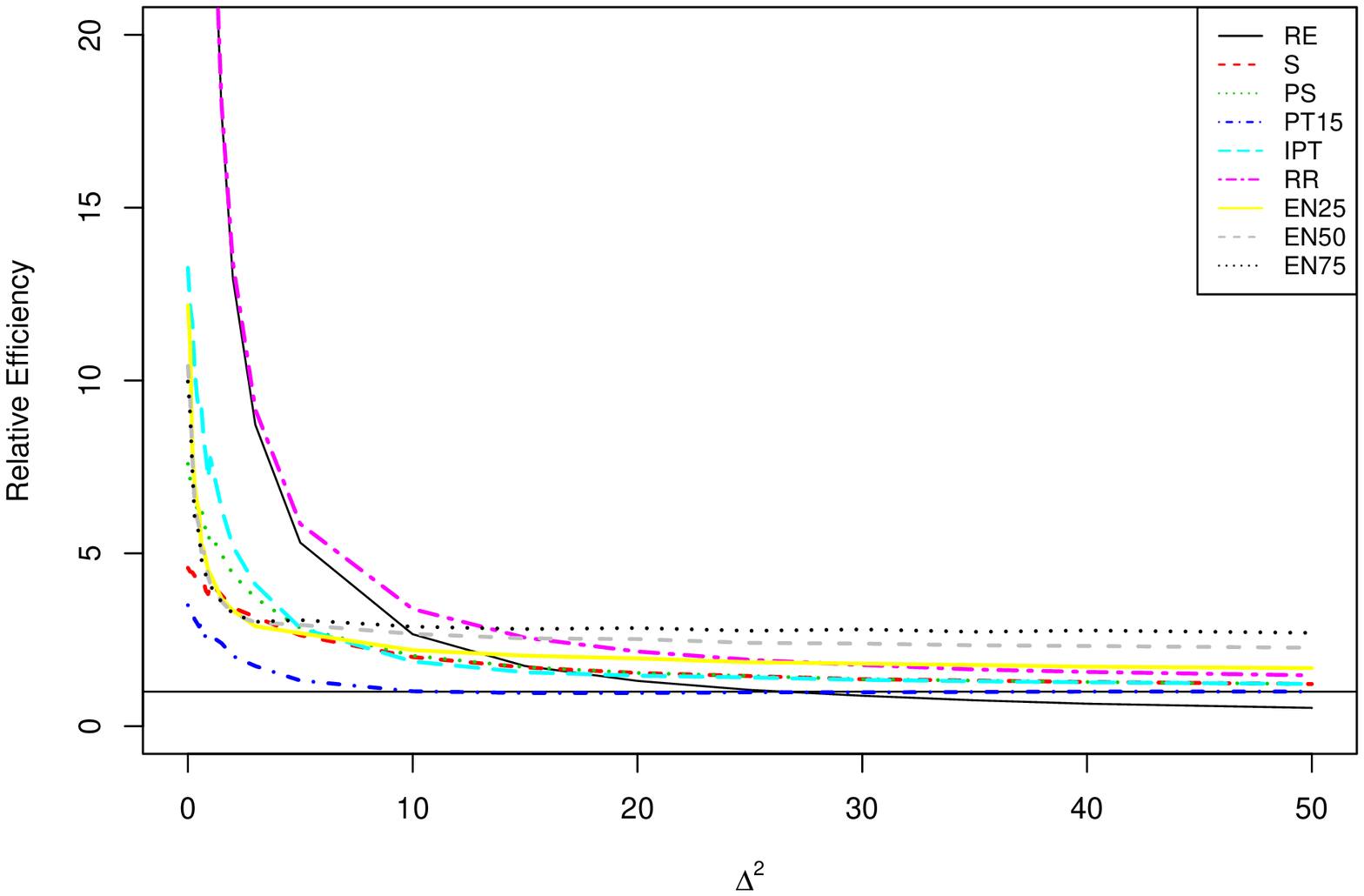}
\caption{Relative Efficiency plots for Shrinkage, Pretest, Ridge, and Elastic Net Estimators when $p=10$, $r=0.2$, and $n=100$.}
\end{figure}

\begin{figure}[!h]\label{fig:ols-ridge-p10-r9}
\includegraphics[width=6in]{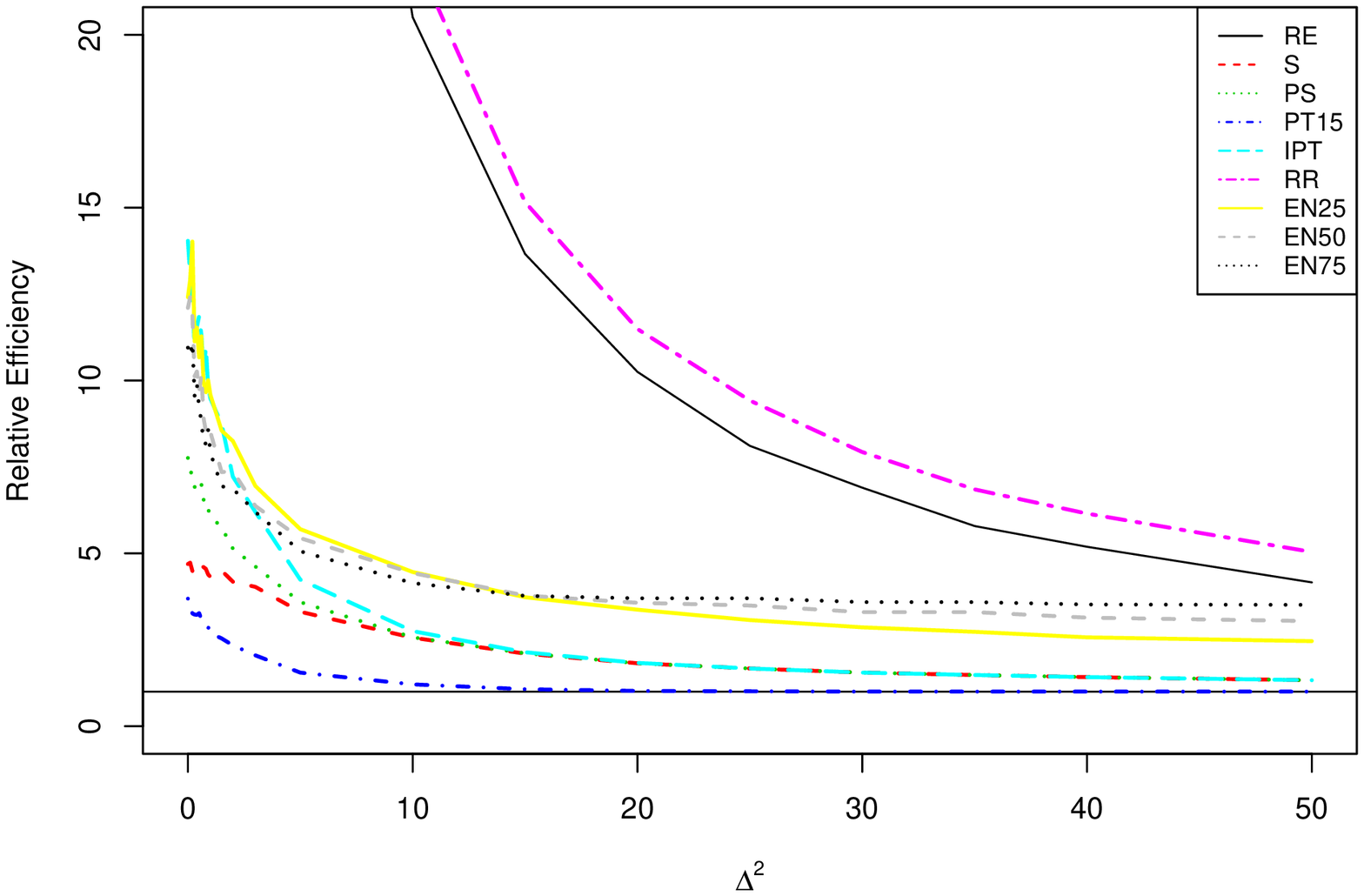}
\caption{Relative Efficiency plots for Shrinkage, Pretest, Ridge, and Elastic Net Estimators when $p=10$, $r=0.9$, and $n=100$.}
\end{figure}

\begin{figure}[!h]
\centering
\includegraphics[width=6in]{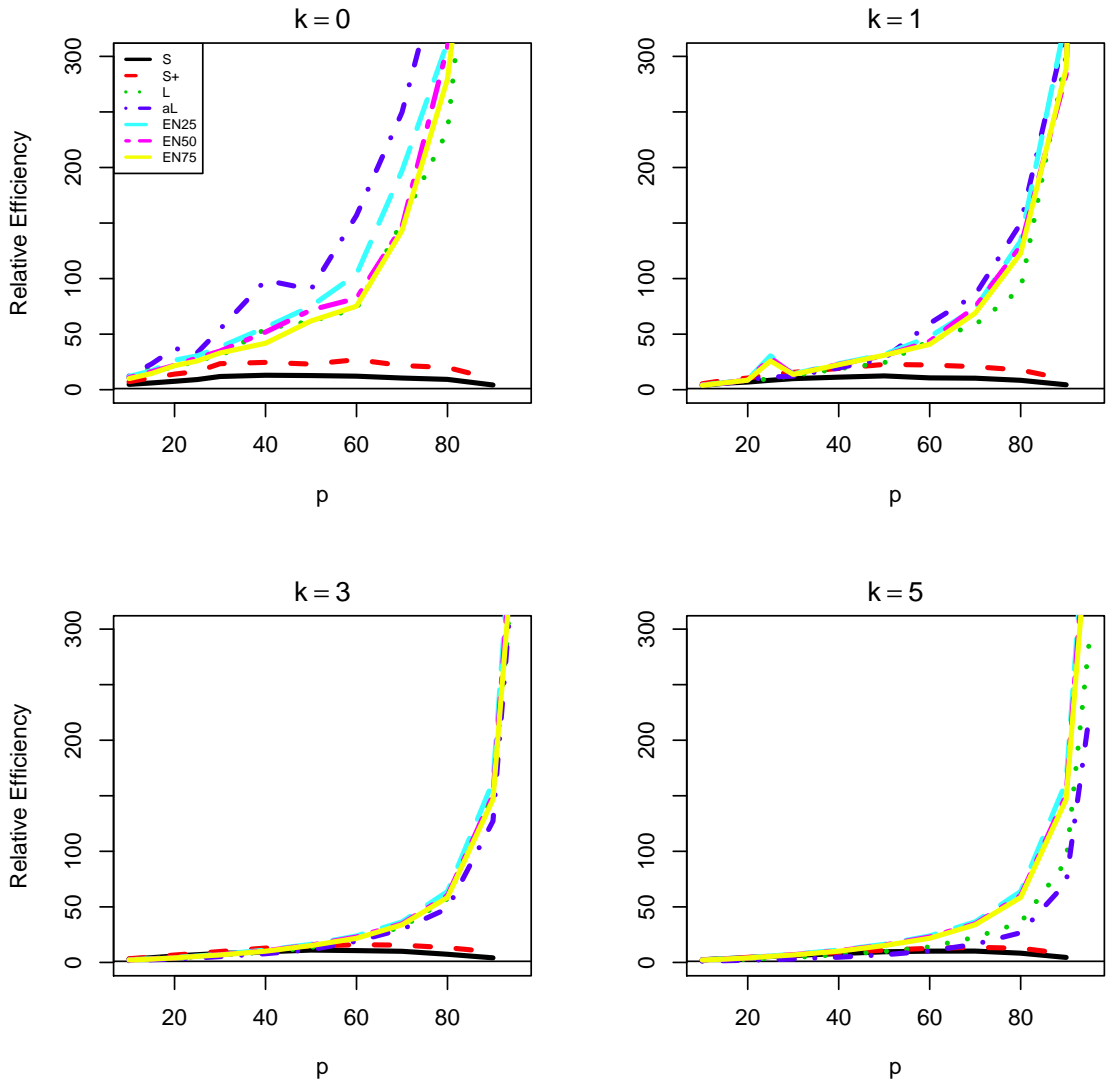}
\caption{Relative Efficiencies of LASSO, aLASSO, ElasticNet, and Stein-type Estimators when $p$ varies and $r=0.2$ for $k=0, 1, 3, 5$. } \label{fig:penalty:JS}
\end{figure}

In this paper, we studied the performance of some selected penalty estimators versus preliminary test and Stein-type estimators under full model hypothesis via simulation studies when the parameter space is smaller than the sample space. We have provided Tables~\ref{tb:rel:eff:cor0}-\ref{tb:penalty:JS:comparision} of relative efficiencies compared to LSE as well as graphical s (\ref{fig:ols-ridge-p10-r0}-\ref{fig:penalty:JS}).

Based on the analysis of the relative efficiencies, we found that (1) ridge regression estimator uniformly dominates the LSE, RE, PTE, JSE and PRSE. (2) LASSO, aLASSO, and SCAD estimators and Stein-type estimators uniformly dominate only LSE. (3) Neither LASSO, aLASSO and SCAD nor the PTE and Stein-type estimator dominate each other group. (4) LASSO and adaptive LASSO outperform SCAD when there are more non-zero coefficients and the covariates are highly correlated. (5) Elastic net is a compromise between LASSO and ridge regression estimators. Relative efficiency of EN is a decreasing function of $\alpha \in (0,1)$ for some $\Delta^2<\Delta^2_{(\alpha, r)}$ depending on the values of $(\alpha,r)$, and increasing for $\Delta^2>\Delta^2_{(\alpha, r)}$. EN75 dominates both types of Stein-type estimators uniformly while EN25 and EN50 dominates S and S+ over a significant portion of the parameter space.

Finally, we note that LASSO is a regularization technique for simultaneous estimation and variable section whereas Stein-type methodology is focused only on estimation which dominates the LSE uniformly under a quadratic risk. 

\begin{landscape}
\begin{table}[!h]
\centering
\caption{Relative efficiencies of Shrinkage, Pretest, Ridge, and Elastic Net Estimators when $p=10$, $r=0$, and $n=100$.} \label{tb:rel:eff:cor0}
\bigskip
\begin{tabular}{rrrrrrrrrrrrrr}
\hline
 $\Delta^2$ & LSE   & RE & \multicolumn{4}{c}{Preliminary Tests} & IPT &  S &  S+ &  RR & EN25 & EN50 & EN75 \\
 &&& $\alpha=.05$ & $\alpha=.15$ & $\alpha=.20$ & $\alpha=.25$ &&  &  & &&& \\
 \hline
   0.0&   1&     Inf&  8.72&  3.69&  2.91&  2.45&  13.78&  4.78&  7.75&      NA& 12.25& 10.64& 10.60 \\
   0.1&   1&  223.95&  7.11&  3.35&  2.82&  2.41&  12.07&  4.63&  7.29&  223.91&  9.40&  8.20&  8.39 \\
   0.2&   1&  113.57&  6.26&  3.18&  2.64&  2.29&  10.84&  4.27&  6.88&  113.74&  7.42&  6.76&  6.60 \\
   0.3&   1&   75.23&  6.84&  3.24&  2.67&  2.25&  11.02&  4.50&  6.89&   75.38&  6.34&  5.81&  5.66 \\
   0.4&   1&   55.98&  6.00&  3.08&  2.54&  2.20&   9.55&  4.35&  6.42&   56.13&  6.23&  5.81&  5.66 \\
   0.5&   1&   45.11&  5.91&  3.07&  2.51&  2.17&   9.55&  4.18&  6.39&   45.18&  5.37&  5.15&  4.95 \\
   0.6&   1&   37.83&  5.70&  2.88&  2.41&  2.08&   8.97&  4.28&  6.13&   38.00&  5.18&  4.82&  4.70 \\
   0.7&   1&   31.42&  5.29&  2.92&  2.49&  2.11&   8.65&  3.99&  6.07&   31.71&  4.65&  4.36&  4.21 \\
   0.8&   1&   28.51&  4.68&  2.66&  2.21&  1.97&   7.69&  3.93&  5.56&   28.54&  4.43&  4.27&  4.07 \\
   0.9&   1&   25.71&  4.38&  2.40&  2.13&  1.91&   7.13&  3.99&  5.34&   25.74&  4.18&  4.02&  3.88 \\
   1.0&   1&   23.12&  4.32&  2.40&  2.06&  1.82&   6.79&  3.86&  5.18&   23.20&  4.04&  3.86&  3.77 \\
   1.5&   1&   15.15&  3.68&  2.27&  1.96&  1.72&   5.98&  3.99&  4.82&   15.48&  3.24&  3.17&  3.12 \\
   2.0&   1&   11.36&  2.99&  1.94&  1.72&  1.53&   4.79&  3.32&  4.15&   11.75&  3.05&  3.04&  3.06 \\
   3.0&   1&    7.52&  2.38&  1.67&  1.49&  1.39&   3.85&  3.09&  3.59&    8.04&  2.75&  2.85&  2.86 \\
   5.0&   1&    4.53&  1.61&  1.27&  1.20&  1.15&   2.63&  2.49&  2.71&    5.04&  2.42&  2.68&  2.79 \\
  10.0&   1&    2.30&  1.04&  0.99&  0.98&  0.98&   1.77&  1.91&  1.96&    2.96&  2.08&  2.47&  2.68 \\
  15.0&   1&    1.51&  0.90&  0.94&  0.95&  0.96&   1.50&  1.64&  1.66&    2.29&  1.94&  2.39&  2.64 \\
  20.0&   1&    1.12&  0.89&  0.95&  0.97&  0.98&   1.38&  1.49&  1.49&    1.96&  1.84&  2.27&  2.61 \\
  25.0&   1&    0.89&  0.91&  0.96&  0.98&  0.98&   1.34&  1.40&  1.40&    1.75&  1.78&  2.28&  2.62 \\
  30.0&   1&    0.76&  0.95&  0.98&  0.99&  0.99&   1.31&  1.33&  1.33&    1.63&  1.71&  2.20&  2.62 \\
  35.0&   1&    0.66&  0.98&  1.00&  1.00&  1.00&   1.28&  1.29&  1.29&    1.56&  1.68&  2.21&  2.59 \\
  40.0&   1&    0.56&  0.99&  1.00&  1.00&  1.00&   1.25&  1.26&  1.26&    1.46&  1.65&  2.20&  2.60 \\
  50.0&   1&    0.45&  1.00&  1.00&  1.00&  1.00&   1.20&  1.20&  1.20&    1.38&  1.59&  2.12&  2.54 \\
\hline
\end{tabular}
\end{table}
\end{landscape}

\begin{landscape}
\begin{table}[!h]
\centering
\caption{Relative efficiencies of Shrinkage, Pretest, Ridge, and Elastic Net Estimators when $p=10$, $r=0.2$, and $n=100$.} \label{tb:rel:eff:cor2}
\bigskip
\begin{tabular}{rrrrrrrrrrrrrr}
\hline
 $\Delta^2$ & LSE   & RE & \multicolumn{4}{c}{Preliminary Tests} & IPT &  S &  S+ &  RR & EN25 & EN50 & EN75 \\
 &&& $\alpha=.05$ & $\alpha=.15$ & $\alpha=.20$ & $\alpha=.25$ &&  &  &&&& \\
 \hline
   0.0&   1&    Inf&  7.69&  3.50&  2.87&  2.42&  13.26&  4.58&  7.59&     Inf & 12.17& 10.42& 9.97\\
   0.1&   1& 259.78&  7.15&  3.32&  2.69&  2.24&  12.14&  4.47&  7.15&  260.34 & 11.06&  9.47& 9.06\\
   0.2&   1& 130.36&  7.42&  3.27&  2.72&  2.30&  11.67&  4.45&  7.03&  130.47 &  8.00&  7.77& 7.38\\
   0.3&   1&  88.84&  6.11&  3.09&  2.54&  2.21&  10.33&  4.35&  6.64&   88.60 &  6.99&  6.49& 6.04\\
   0.4&   1&  66.38&  5.69&  2.99&  2.48&  2.08&   9.56&  4.23&  6.33&   66.57 &  6.50&  6.03& 5.82\\
   0.5&   1&  52.90&  5.55&  2.88&  2.41&  2.13&   9.14&  4.24&  6.22&   53.12 &  6.21&  5.77& 5.55\\
   0.6&   1&  42.85&  5.90&  2.96&  2.49&  2.18&   9.24&  4.12&  6.34&   43.07 &  5.27&  5.04& 4.84\\
   0.7&   1&  37.44&  5.23&  2.66&  2.28&  2.04&   8.32&  4.17&  5.91&   37.58 &  5.13&  5.05& 4.56\\
   0.8&   1&  33.14&  4.71&  2.58&  2.21&  1.96&   7.82&  3.92&  5.67&   33.32 &  4.82&  4.72& 4.37\\
   0.9&   1&  29.40&  4.39&  2.53&  2.14&  1.89&   7.17&  3.82&  5.39&   29.69 &  4.52&  4.32& 4.27\\
   1.0&   1&  26.29&  4.87&  2.62&  2.19&  1.90&   7.76&  4.15&  5.55&   26.66 &  4.40&  4.11& 4.06\\
   1.5&   1&  17.75&  3.83&  2.40&  2.07&  1.82&   6.29&  3.79&  4.97&   17.96 &  3.69&  3.61& 3.59\\
   2.0&   1&  12.96&  3.22&  2.03&  1.80&  1.63&   5.20&  3.43&  4.40&   13.43 &  3.35&  3.26& 3.26\\
   3.0&   1&   8.72&  2.48&  1.74&  1.55&  1.44&   4.10&  3.16&  3.71&    9.16 &  2.89&  3.00& 3.02\\
   5.0&   1&   5.31&  1.68&  1.32&  1.23&  1.18&   2.84&  2.63&  2.84&    5.85 &  2.69&  2.93& 3.07\\
  10.0&   1&   2.66&  1.08&  1.01&  1.00&  1.00&   1.87&  2.00&  2.04&    3.39 &  2.20&  2.67& 2.88\\
  15.0&   1&   1.74&  0.93&  0.96&  0.96&  0.97&   1.56&  1.70&  1.71&    2.56 &  2.04&  2.54& 2.81\\
  20.0&   1&   1.31&  0.92&  0.96&  0.97&  0.98&   1.46&  1.54&  1.54&    2.16 &  1.96&  2.52& 2.84\\
  25.0&   1&   1.05&  0.94&  0.98&  0.99&  0.99&   1.41&  1.45&  1.45&    1.92 &  1.85&  2.41& 2.76\\
  30.0&   1&   0.88&  0.95&  0.98&  0.99&  0.99&   1.34&  1.36&  1.36&    1.77 &  1.81&  2.39& 2.80\\
  35.0&   1&   0.75&  0.97&  0.99&  1.00&  1.00&   1.30&  1.32&  1.32&    1.64 &  1.77&  2.34& 2.73\\
  40.0&   1&   0.65&  0.99&  1.00&  1.00&  1.00&   1.27&  1.28&  1.28&    1.57 &  1.72&  2.32& 2.77\\
  50.0&   1&   0.53&  1.00&  1.00&  1.00&  1.00&   1.22&  1.22&  1.22&    1.47 &  1.68&  2.27& 2.70\\
  \hline
\end{tabular}
\end{table}
\end{landscape}

\begin{landscape}
\begin{table}[!h]
\centering
\caption{Relative efficiencies of Shrinkage, Pretest, Ridge, and Elastic Net Estimators when $p=10$, $r=0.9$, and $n=100$..} \label{tb:rel:eff:cor9}
\bigskip
\begin{tabular}{rrrrrrrrrrrrrr}
\hline
 $\Delta^2$ & LSE   & RE & \multicolumn{4}{c}{Preliminary Tests} & IPT &  S &  S+ &  RR & EN25 & EN50 & EN75 \\
 &&& $\alpha=.05$ & $\alpha=.15$ & $\alpha=.20$ & $\alpha=.25$ &&  &  &&&& \\
 \hline
   0.0&   1&      Inf& 8.38& 3.69& 2.93& 2.44& 14.04& 4.69& 7.76&     Inf& 12.41& 12.10& 10.95 \\
   0.1&   1&  2039.90& 7.44& 3.48& 2.75& 2.35& 13.06& 4.73& 7.43& 2040.14& 12.91& 12.41& 10.72 \\
   0.2&   1&  1016.40& 7.51& 3.26& 2.68& 2.31& 12.31& 4.49& 7.22& 1023.32& 14.02& 11.88& 11.10 \\
   0.3&   1&   684.67& 6.28& 3.23& 2.66& 2.29& 11.00& 4.53& 6.92&  687.72& 11.13& 10.12&  9.55 \\
   0.4&   1&   515.71& 6.69& 3.22& 2.67& 2.29& 11.55& 4.44& 6.94&  518.54& 11.53& 10.26&  9.90 \\
   0.5&   1&   396.49& 7.45& 3.27& 2.66& 2.28& 11.86& 4.43& 7.13&  401.43& 10.67&  9.66&  9.22 \\
   0.6&   1&   331.07& 6.65& 3.10& 2.66& 2.30& 11.32& 4.64& 7.00&  335.24& 11.28& 10.19&  8.82 \\
   0.7&   1&   293.18& 5.73& 2.99& 2.47& 2.13& 10.28& 4.60& 6.55&  298.05&  9.93&  8.90&  8.17 \\
   0.8&   1&   256.39& 6.47& 2.97& 2.46& 2.12& 10.97& 4.55& 6.58&  260.43&  9.66&  8.53&  8.18 \\
   0.9&   1&   227.77& 5.67& 2.87& 2.33& 2.07&  9.84& 4.39& 6.28&  232.61& 10.01&  8.71&  8.58 \\
   1.0&   1&   202.78& 5.54& 2.69& 2.28& 1.99&  9.48& 4.32& 6.08&  207.87&  9.59&  8.51&  7.91 \\
   1.5&   1&   137.01& 5.11& 2.54& 2.19& 1.92&  8.74& 4.50& 5.75&  141.35&  8.55&  7.35&  6.94 \\
   2.0&   1&   103.02& 4.12& 2.34& 1.97& 1.74&  7.22& 4.18& 5.14&  107.37&  8.25&  7.38&  6.87 \\
   3.0&   1&    67.99& 3.40& 2.05& 1.77& 1.60&  6.20& 4.03& 4.62&   71.81&  6.95&  6.37&  6.20 \\
   5.0&   1&    40.70& 2.31& 1.55& 1.41& 1.33&  4.24& 3.31& 3.59&   43.66&  5.70&  5.44&  5.06 \\
  10.0&   1&    20.51& 1.47& 1.21& 1.15& 1.12&  2.75& 2.56& 2.57&   22.41&  4.46&  4.43&  4.14 \\
  15.0&   1&    13.66& 1.19& 1.07& 1.06& 1.04&  2.14& 2.10& 2.11&   15.16&  3.73&  3.78&  3.77 \\
  20.0&   1&    10.25& 1.07& 1.02& 1.02& 1.01&  1.83& 1.82& 1.82&   11.49&  3.37&  3.57&  3.70 \\
  25.0&   1&     8.11& 1.02& 1.01& 1.00& 1.00&  1.67& 1.67& 1.67&    9.42&  3.07&  3.49&  3.70 \\
  30.0&   1&     6.90& 1.01& 1.00& 1.00& 1.00&  1.55& 1.55& 1.55&    7.93&  2.86&  3.30&  3.59 \\
  35.0&   1&     5.79& 1.01& 1.00& 1.00& 1.00&  1.48& 1.48& 1.48&    6.85&  2.73&  3.30&  3.59 \\
  40.0&   1&     5.19& 1.00& 1.00& 1.00& 1.00&  1.42& 1.42& 1.42&    6.15&  2.57&  3.14&  3.52 \\
  50.0&   1&     4.16& 1.00& 1.00& 1.00& 1.00&  1.33& 1.33& 1.33&    5.04&  2.46&  3.04&  3.51 \\
  \hline
\end{tabular}
\end{table}
\end{landscape}

\begin{landscape}
\begin{table}[!h]
\centering
\caption{Relative efficiencies of LASSO, aLASSO, and SCAD estimators when $n=100$, $r=0.2,$ $p=10$.} \label{tb:rel:eff:penalty:cor2:p10}
\bigskip
\begin{tabular}{r|rrr|rrr|rrr|rrr|rrr}
\hline
 $\Delta^2$ & \multicolumn{3}{c|}{$k=1$} & \multicolumn{3}{c|}{$k=2$} & \multicolumn{3}{c|}{$k=3$}
 & \multicolumn{3}{c|}{$k=4$} & \multicolumn{3}{c}{$k=5$} \\
 & L & aL & SCAD & L & aL & SCAD & L & aL & SCAD & L & aL & SCAD & L & aL & SCAD \\
 \hline
   0.0& 4.16& 4.46& 3.27& 2.56& 2.05& 1.83&2.07& 1.48& 1.36&1.74& 1.25& 1.12&1.49& 1.09& 0.96\\  
   0.1& 3.52& 3.52& 2.50& 2.32& 1.84& 1.67&1.95& 1.51& 1.44&1.71& 1.24& 1.11&1.43& 1.11& 0.97\\  
   0.2& 3.52& 3.43& 3.19& 2.47& 1.81& 1.68&1.88& 1.39& 1.37&1.65& 1.23& 1.10&1.46& 1.08& 0.95\\  
   0.3& 3.03& 2.96& 2.90& 2.32& 2.01& 1.80&1.79& 1.32& 1.27&1.59& 1.11& 1.06&1.49& 1.10& 1.00\\  
   0.4& 3.30& 2.88& 2.56& 2.19& 1.86& 1.56&1.71& 1.38& 1.25&1.62& 1.18& 1.09&1.40& 1.06& 0.92\\  
   0.5& 2.72& 2.37& 2.35& 2.23& 1.73& 1.62&1.79& 1.27& 1.22&1.53& 1.15& 1.02&1.34& 1.04& 0.93\\  
   0.6& 2.90& 2.43& 2.22& 2.29& 1.73& 1.60&1.70& 1.43& 1.18&1.56& 1.13& 1.03&1.44& 1.05& 0.95\\  
   0.7& 2.91& 2.39& 2.08& 2.14& 1.58& 1.48&1.75& 1.28& 1.19&1.53& 1.16& 0.99&1.37& 1.08& 0.94\\  
   0.8& 2.68& 2.43& 2.23& 2.21& 1.56& 1.47&1.70& 1.29& 1.13&1.46& 1.11& 1.02&1.30& 1.03& 0.90\\  
   0.9& 2.51& 2.04& 1.83& 1.89& 1.69& 1.44&1.70& 1.22& 1.16&1.58& 1.23& 1.07&1.39& 1.04& 0.92\\  
   1.0& 2.85& 2.30& 1.97& 1.88& 1.51& 1.29&1.68& 1.27& 1.16&1.49& 1.19& 0.99&1.37& 1.05& 0.88\\  
   1.5& 2.30& 1.91& 1.75& 2.00& 1.48& 1.29&1.76& 1.26& 1.14&1.51& 1.20& 0.99&1.38& 1.11& 0.93\\  
   2.0& 2.37& 1.95& 1.94& 2.03& 1.60& 1.37&1.78& 1.27& 1.14&1.55& 1.15& 1.00&1.32& 1.02& 0.90\\  
   3.0& 2.16& 2.00& 1.81& 1.77& 1.43& 1.26&1.63& 1.34& 1.11&1.47& 1.17& 1.00&1.33& 1.08& 0.94\\  
   5.0& 2.22& 2.10& 1.87& 1.87& 1.67& 1.37&1.58& 1.35& 1.20&1.53& 1.21& 1.02&1.28& 1.04& 0.93\\  
  10.0& 2.31& 2.34& 2.03& 2.02& 1.68& 1.47&1.57& 1.31& 1.23&1.40& 1.13& 1.05&1.30& 1.01& 0.92\\  
  15.0& 2.10& 2.24& 2.33& 1.91& 1.71& 1.62&1.73& 1.38& 1.14&1.50& 1.19& 1.07&1.34& 1.03& 0.91\\  
  20.0& 2.22& 2.37& 2.61& 1.81& 1.61& 1.56&1.50& 1.26& 1.24&1.45& 1.20& 1.08&1.29& 0.99& 0.91\\  
  25.0& 2.39& 2.67& 2.67& 2.04& 1.81& 1.69&1.67& 1.42& 1.30&1.50& 1.20& 1.06&1.32& 1.06& 0.98\\  
  30.0& 2.32& 2.56& 2.46& 1.80& 1.65& 1.59&1.52& 1.25& 1.24&1.45& 1.14& 0.99&1.34& 1.01& 0.90\\  
  35.0& 2.36& 2.71& 2.81& 1.87& 1.72& 1.67&1.69& 1.37& 1.25&1.43& 1.13& 1.02&1.33& 1.03& 0.90\\  
  40.0& 2.30& 2.62& 2.37& 2.00& 1.88& 1.75&1.63& 1.43& 1.32&1.43& 1.07& 0.99&1.32& 1.11& 0.95\\  
  50.0& 2.30& 2.42& 2.35& 2.10& 1.98& 1.93&1.62& 1.36& 1.27&1.45& 1.10& 1.04&1.31& 1.02& 0.89\\
  \hline
\end{tabular}
\end{table}
\end{landscape}

\begin{landscape}
\begin{table}[!h]
\centering
\caption{Relative efficiencies of LASSO, aLASSO, and SCAD estimators when $n=100$, $r=0.9,$ $p=10$.} \label{tb:rel:eff:penalty:cor9:p10}
\bigskip
\begin{tabular}{r|rrr|rrr|rrr|rrr|rrr}
\hline
 $\Delta^2$ & \multicolumn{3}{c|}{$k=1$} & \multicolumn{3}{c|}{$k=2$} & \multicolumn{3}{c|}{$k=3$}
 & \multicolumn{3}{c|}{$k=4$} & \multicolumn{3}{c}{$k=5$} \\
 & L & aL & SCAD & L & aL & SCAD & L & aL & SCAD & L & aL & SCAD & L & aL & SCAD \\
 \hline
   0.0& 9.14& 12.00& 3.56&5.10& 5.13& 2.49&4.51& 3.66& 2.30&4.05& 3.38& 1.72&3.51& 2.20& 1.36 \\  
   0.1& 6.57&  8.72& 2.78&5.28& 5.09& 2.37&4.25& 3.37& 2.17&3.81& 2.68& 1.50&3.11& 2.18& 1.33 \\  
   0.2& 8.04&  8.76& 3.52&5.03& 4.95& 2.53&4.16& 3.03& 1.80&3.47& 2.74& 1.55&3.36& 2.07& 1.28 \\  
   0.3& 6.47&  7.66& 2.71&5.96& 4.85& 2.36&3.74& 3.46& 1.68&3.11& 2.37& 1.47&2.84& 2.13& 1.27 \\  
   0.4& 8.04& 11.11& 3.50&5.11& 4.23& 2.22&4.33& 3.18& 1.76&3.56& 2.61& 1.35&2.95& 2.02& 1.21 \\  
   0.5& 6.25&  6.45& 2.69&5.28& 4.50& 2.38&3.47& 2.83& 1.67&3.58& 2.34& 1.39&2.76& 1.89& 1.21 \\  
   0.6& 6.64&  6.41& 3.16&5.06& 4.00& 2.16&4.22& 3.02& 1.79&3.74& 2.69& 1.47&2.95& 2.04& 1.24 \\  
   0.7& 5.58&  6.29& 2.70&6.33& 5.12& 2.53&3.49& 2.54& 1.51&3.34& 2.23& 1.32&2.95& 2.07& 1.23 \\  
   0.8& 5.31&  7.69& 2.69&4.73& 4.22& 2.04&3.48& 2.46& 1.51&3.45& 2.23& 1.27&3.09& 2.10& 1.27 \\  
   0.9& 5.04&  4.90& 2.70&5.78& 4.38& 1.87&4.32& 2.65& 1.64&3.14& 2.30& 1.33&2.71& 1.85& 1.17 \\  
   1.0& 5.61&  6.59& 2.66&5.59& 4.05& 2.17&3.69& 3.17& 1.42&3.46& 2.73& 1.49&2.71& 1.86& 1.15 \\  
   1.5& 5.18&  4.64& 2.56&4.70& 3.38& 2.08&3.57& 2.63& 1.46&3.40& 2.41& 1.30&2.64& 1.93& 1.10 \\  
   2.0& 5.06&  4.77& 2.35&3.91& 3.67& 1.66&4.01& 2.73& 1.64&2.94& 2.08& 1.20&2.66& 1.80& 1.14 \\  
   3.0& 4.69&  4.96& 2.51&3.83& 2.79& 1.72&3.17& 2.27& 1.39&2.85& 2.04& 1.24&2.64& 1.80& 1.12 \\  
   5.0& 4.51&  3.55& 1.96&2.92& 2.63& 1.42&2.95& 2.17& 1.39&2.82& 1.95& 1.20&2.35& 1.55& 1.02 \\  
  10.0& 3.70&  3.71& 1.69&2.93& 2.46& 1.56&2.97& 2.24& 1.43&2.81& 2.06& 1.23&2.41& 1.62& 1.07 \\  
  15.0& 3.47&  3.47& 1.67&2.61& 2.20& 1.49&3.17& 2.27& 1.38&2.70& 1.81& 1.17&2.28& 1.81& 1.20 \\  
  20.0& 3.16&  3.27& 1.99&3.41& 2.80& 1.76&2.91& 2.32& 1.66&2.39& 1.66& 1.22&2.41& 1.66& 1.20 \\  
  25.0& 3.40&  3.31& 2.07&3.04& 2.70& 1.83&2.94& 2.17& 1.79&2.36& 2.03& 1.41&2.21& 1.72& 1.17 \\  
  30.0& 3.19&  4.62& 2.22&3.17& 3.22& 1.96&2.81& 2.44& 1.82&2.73& 2.21& 1.47&2.30& 1.64& 1.29 \\  
  35.0& 3.52&  5.22& 2.56&2.88& 2.79& 1.97&2.59& 2.35& 1.93&2.46& 1.93& 1.37&2.45& 1.72& 1.28 \\  
  40.0& 3.39&  6.06& 2.37&3.33& 3.89& 3.14&2.55& 2.31& 1.80&2.46& 1.91& 1.50&2.33& 1.64& 1.33 \\  
  50.0& 3.84&  5.48& 3.90&3.39& 3.50& 2.73&3.04& 2.60& 1.85&2.49& 2.04& 1.43&2.28& 1.71& 1.26 \\
  \hline
\end{tabular}
\end{table}
\end{landscape}

\begin{landscape}
\begin{table}[!h]
\centering
\caption{Relative efficiencies of LASSO, aLASSO, and SCAD estimators when $n=100$, $r=0.2,$ $p=20$.} \label{tb:rel:eff:penalty:p20cor2}
\bigskip
\begin{tabular}{r|rrr|rrr|rrr|rrr|rrr}
\hline
 $\Delta^2$ & \multicolumn{3}{c|}{$k=1$} & \multicolumn{3}{c|}{$k=2$} & \multicolumn{3}{c|}{$k=3$}
 & \multicolumn{3}{c|}{$k=4$} & \multicolumn{3}{c}{$k=5$} \\
 & L & aL & SCAD & L & aL & SCAD & L & aL & SCAD & L & aL & SCAD & L & aL & SCAD \\
 \hline
   0.0& 10.04& 9.48& 11.55& 5.89& 5.32& 5.21& 3.95& 3.01& 3.42& 3.13& 2.24& 2.60& 2.83& 2.02& 2.03 \\  
   0.1&  8.23& 6.94&  8.00& 5.07& 4.57& 4.50& 4.01& 3.14& 3.60& 3.29& 2.28& 2.60& 2.81& 1.94& 1.72 \\  
   0.2&  6.78& 6.42&  8.41& 5.18& 4.02& 3.97& 3.64& 2.89& 2.99& 3.37& 2.31& 2.55& 2.79& 1.87& 1.70 \\  
   0.3&  7.77& 8.07&  7.67& 5.05& 3.90& 4.73& 3.65& 2.70& 2.54& 3.14& 2.19& 2.17& 2.65& 1.85& 1.78 \\  
   0.4&  6.98& 6.15&  7.49& 4.46& 3.61& 3.98& 3.53& 2.54& 3.00& 3.07& 2.22& 2.17& 2.78& 1.98& 1.88 \\  
   0.5&  5.79& 5.12&  6.04& 4.71& 3.90& 3.61& 3.48& 2.63& 2.98& 2.99& 2.08& 2.05& 2.57& 1.85& 1.84 \\  
   0.6&  5.71& 5.83&  6.12& 4.28& 3.51& 4.11& 3.40& 2.33& 2.63& 3.06& 2.27& 2.28& 2.55& 1.72& 1.52 \\  
   0.7&  4.97& 5.32&  4.54& 4.18& 3.38& 3.54& 3.83& 2.53& 2.92& 3.08& 2.15& 2.29& 2.63& 1.88& 1.83 \\  
   0.8&  5.15& 4.77&  4.63& 4.04& 3.45& 3.88& 3.19& 2.40& 2.58& 2.83& 1.97& 2.03& 2.45& 1.66& 1.67 \\  
   0.9&  5.67& 4.71&  5.04& 4.25& 3.38& 3.41& 3.28& 2.42& 2.33& 2.98& 2.10& 1.93& 2.63& 1.82& 1.77 \\  
   1.0&  5.84& 5.52&  5.12& 4.18& 3.06& 3.49& 3.33& 2.47& 2.73& 2.99& 2.10& 2.21& 2.42& 1.68& 1.67 \\  
   1.5&  5.52& 4.55&  5.22& 3.67& 2.96& 3.09& 3.06& 2.25& 2.32& 2.59& 1.92& 1.84& 2.42& 1.67& 1.64 \\  
   2.0&  4.66& 3.86&  4.22& 3.77& 2.67& 2.85& 3.15& 2.35& 2.45& 2.84& 2.04& 2.01& 2.55& 1.77& 1.58 \\  
   3.0&  4.69& 4.08&  4.35& 3.60& 2.86& 2.87& 3.09& 2.26& 2.56& 2.64& 2.05& 1.99& 2.44& 1.74& 1.52 \\  
   5.0&  4.26& 3.79&  4.11& 3.83& 2.81& 2.86& 3.04& 2.48& 2.28& 2.58& 1.97& 1.83& 2.44& 1.69& 1.57 \\  
  10.0&  4.66& 6.26&  4.97& 3.54& 3.08& 3.77& 3.21& 2.61& 2.85& 2.80& 2.23& 2.15& 2.38& 1.84& 1.74 \\  
  15.0&  4.15& 5.86&  6.78& 3.76& 4.48& 3.41& 2.94& 2.54& 2.96& 2.74& 2.26& 2.27& 2.50& 1.87& 1.80 \\  
  20.0&  4.01& 5.44&  6.09& 3.83& 4.03& 4.17& 2.92& 2.48& 2.84& 2.88& 2.43& 2.59& 2.14& 1.68& 1.61 \\  
  25.0&  4.93& 5.57&  6.43& 3.61& 3.85& 4.17& 2.88& 2.37& 2.75& 2.76& 2.28& 2.40& 2.38& 1.84& 1.94 \\  
  30.0&  4.58& 5.78&  6.55& 3.64& 3.54& 4.14& 3.03& 2.89& 3.04& 2.73& 2.07& 2.51& 2.43& 1.85& 1.69 \\  
  35.0&  4.45& 5.95&  5.87& 3.41& 3.78& 4.77& 3.34& 3.23& 3.26& 2.67& 2.13& 2.10& 2.41& 1.85& 1.77 \\  
  40.0&  4.34& 5.23&  7.90& 3.66& 3.66& 3.82& 2.95& 2.91& 3.09& 2.71& 2.24& 2.44& 2.33& 1.77& 1.74 \\  
  50.0&  4.34& 6.34&  5.54& 3.72& 3.70& 3.70& 3.12& 2.47& 3.07& 2.75& 2.17& 2.29& 2.42& 1.86& 1.89 \\
  \hline
\end{tabular}
\end{table}
\end{landscape}

\begin{landscape}
\begin{table}[!h]
\centering
\caption{Relative efficiencies of LASSO, aLASSO, and SCAD estimators when $n=100$, $r=0.9,$ $p=20$.} \label{tb:rel:eff:penalty:p20cor9}
\bigskip
\begin{tabular}{r|rrr|rrr|rrr|rrr|rrr}
\hline
 $\Delta^2$ & \multicolumn{3}{c|}{$k=1$} & \multicolumn{3}{c|}{$k=2$} & \multicolumn{3}{c|}{$k=3$}
 & \multicolumn{3}{c|}{$k=4$} & \multicolumn{3}{c}{$k=5$} \\
 & L & aL & SCAD & L & aL & SCAD & L & aL & SCAD & L & aL & SCAD & L & aL & SCAD \\
 \hline
   0.0& 23.92&  24.02&   8.26&  12.45&  16.39&  5.32&  11.14&  7.79&  4.60&  7.54&  6.12&  3.21&  6.91&  4.74&  2.68 \\  
   0.1& 20.15&  27.39&  10.85&  11.67&   8.97&  4.96&   7.55&  6.10&  3.66&  6.95&  5.37&  2.68&  6.63&  4.64&  2.42 \\  
   0.2& 12.21&  12.23&   9.62&   9.57&   7.97&  3.81&   8.71&  6.78&  3.51&  7.98&  4.87&  2.99&  6.23&  3.80&  2.18 \\  
   0.3& 11.47&  14.23&   8.56&   8.89&   7.76&  5.34&   7.89&  6.65&  3.51&  7.20&  5.53&  2.79&  6.03&  4.06&  2.19 \\  
   0.4& 15.30&  14.01&   7.45&   9.61&   9.82&  4.67&   7.74&  5.56&  3.35&  7.81&  4.82&  2.63&  6.12&  4.08&  2.16 \\  
   0.5& 14.67&  11.81&   8.92&   8.05&   7.72&  4.22&   8.40&  5.80&  3.40&  8.36&  4.84&  2.86&  6.73&  4.34&  2.11 \\  
   0.6& 10.71&  13.74&   5.22&  10.31&   9.24&  3.97&   9.55&  6.18&  3.85&  6.64&  4.55&  2.37&  6.22&  3.82&  2.19 \\  
   0.7&  8.82&   9.75&   4.33&   9.53&   8.81&  4.81&   7.82&  6.63&  3.95&  7.23&  4.67&  2.65&  6.61&  4.10&  2.35 \\  
   0.8& 13.56&  12.46&   8.19&  11.37&  10.04&  4.51&   6.67&  6.01&  3.49&  7.89&  4.96&  2.76&  5.79&  3.89&  2.25 \\  
   0.9& 15.31&  15.58&   9.34&  10.71&   8.26&  4.34&   7.18&  5.23&  2.84&  6.77&  4.37&  2.51&  5.18&  3.93&  2.10 \\  
   1.0& 13.25&  10.62&   6.70&   7.96&   7.49&  4.43&   7.55&  5.28&  2.99&  6.65&  4.79&  2.48&  5.12&  3.57&  2.16 \\  
   1.5& 11.19&  10.46&   5.30&   7.15&   5.99&  3.38&   6.93&  5.00&  2.72&  5.75&  4.67&  2.59&  5.33&  3.29&  2.09 \\  
   2.0& 11.93&   9.80&   4.78&   8.88&   6.17&  3.79&   6.63&  5.34&  2.62&  6.20&  4.34&  2.50&  6.01&  3.98&  2.21 \\  
   3.0&  8.50&   7.02&   3.87&   6.34&   4.71&  2.87&   6.04&  4.97&  2.76&  4.95&  3.37&  2.09&  5.02&  3.46&  1.82 \\  
   5.0&  6.78&   5.38&   3.85&   5.99&   5.46&  2.73&   6.00&  3.84&  2.08&  5.84&  4.03&  2.28&  4.95&  3.37&  1.85 \\  
  10.0&  6.25&   4.90&   2.74&   5.31&   4.25&  2.55&   5.77&  4.06&  2.25&  4.99&  3.37&  2.04&  4.11&  3.03&  1.80 \\  
  15.0&  7.55&   7.00&   3.13&   5.45&   4.64&  2.71&   5.44&  4.20&  2.55&  4.72&  3.03&  1.93&  4.42&  3.32&  1.94 \\  
  20.0&  8.01&   7.07&   4.77&   5.59&   4.72&  2.81&   5.25&  4.04&  2.62&  4.80&  3.63&  2.13&  4.63&  3.12&  2.23 \\  
  25.0&  7.23&   6.44&   4.42&   5.46&   4.81&  3.46&   5.49&  4.94&  3.19&  4.57&  3.51&  2.44&  4.06&  2.80&  1.73 \\  
  30.0&  8.61&  10.54&   6.47&   5.50&   5.62&  3.67&   5.28&  4.32&  3.07&  4.41&  3.64&  2.58&  4.58&  3.30&  2.02 \\  
  35.0&  8.15&  12.24&   8.14&   5.81&   5.87&  3.49&   5.29&  4.61&  3.38&  4.36&  3.79&  2.56&  4.12&  3.05&  2.12 \\  
  40.0&  6.64&   9.71&   4.41&   5.83&   6.46&  4.01&   4.80&  4.72&  3.34&  4.44&  3.85&  2.42&  4.34&  3.32&  1.97 \\  
  50.0&  7.67&  13.28&   7.91&   5.61&   6.39&  4.86&   5.76&  5.95&  4.01&  4.68&  3.66&  2.91&  4.16&  3.01&  2.13 \\
  \hline
\end{tabular}
\end{table}
\end{landscape}

\begin{landscape}
\begin{table}[!h]
\centering
\caption{Relative efficiencies of LASSO, aLASSO, and SCAD estimators when $n=100$, $r=0.2,$ $p=30$.} \label{tb:rel:eff:penalty:p30cor2}
\bigskip
\begin{tabular}{r|rrr|rrr|rrr|rrr|rrr}
\hline
 $\Delta^2$ & \multicolumn{3}{c|}{$k=1$} & \multicolumn{3}{c|}{$k=2$} & \multicolumn{3}{c|}{$k=3$}
 & \multicolumn{3}{c|}{$k=4$} & \multicolumn{3}{c}{$k=5$} \\
 & L & aL & SCAD & L & aL & SCAD & L & aL & SCAD & L & aL & SCAD & L & aL & SCAD \\
 \hline
  0.0& 16.80& 13.93& 16.77& 8.69& 7.44&  9.41& 7.61& 5.63& 6.86& 5.08& 4.05& 5.36& 4.16& 3.07& 3.97\\  
  0.1& 11.19&  9.76& 13.00& 8.33& 6.48& 10.39& 5.93& 4.66& 5.96& 5.10& 3.99& 4.67& 4.42& 3.18& 3.86\\  
  0.2& 13.34& 12.73& 16.42& 7.31& 5.92&  7.01& 6.22& 5.30& 5.83& 5.10& 3.65& 4.39& 4.06& 2.88& 3.69\\  
  0.3& 13.21& 12.29& 13.85& 7.40& 7.10&  7.60& 5.80& 4.66& 5.84& 4.91& 3.59& 4.12& 4.18& 2.92& 3.36\\  
  0.4&  9.48& 10.63&  8.38& 7.14& 6.73&  8.08& 5.72& 4.30& 5.17& 4.85& 3.38& 4.67& 4.22& 2.95& 3.60\\  
  0.5& 10.83&  9.44& 12.30& 8.10& 6.67&  8.20& 5.58& 4.55& 5.48& 4.28& 3.07& 3.70& 4.25& 2.79& 2.89\\  
  0.6&  9.43&  8.75& 11.54& 7.00& 5.49&  7.91& 6.09& 4.48& 5.74& 4.34& 3.17& 3.94& 4.00& 2.72& 3.09\\  
  0.7& 10.00&  8.75&  9.88& 6.30& 5.70&  6.06& 5.70& 4.36& 5.74& 4.31& 3.32& 3.87& 3.98& 2.78& 3.09\\  
  0.8&  9.77&  9.11& 11.12& 7.25& 5.63&  7.27& 5.57& 3.95& 5.08& 4.66& 3.27& 4.00& 4.28& 2.72& 3.32\\  
  0.9& 10.13&  9.43&  9.24& 6.59& 5.01&  6.53& 5.10& 3.96& 5.40& 4.53& 3.33& 3.58& 3.74& 2.67& 2.78\\  
  1.0&  8.41&  7.99& 10.96& 6.68& 5.26&  6.85& 5.07& 3.75& 4.34& 4.66& 3.36& 3.47& 3.82& 2.58& 3.09\\  
  1.5&  7.99&  6.25&  8.31& 6.31& 5.09&  5.55& 4.92& 3.60& 4.73& 4.47& 3.08& 3.80& 3.85& 2.61& 3.03\\  
  2.0&  7.37&  6.68&  8.42& 5.83& 4.57&  5.60& 5.20& 4.19& 4.79& 4.33& 3.12& 3.10& 3.63& 2.59& 2.83\\  
  3.0&  7.68&  6.43&  7.81& 5.62& 4.44&  5.30& 4.98& 4.11& 4.46& 4.34& 3.36& 3.86& 3.89& 2.70& 2.91\\  
  5.0&  7.57&  7.24&  7.24& 6.04& 5.49&  5.74& 4.91& 3.96& 4.54& 4.40& 3.33& 3.69& 3.66& 2.78& 2.81\\  
 10.0&  8.11&  9.32& 10.06& 6.31& 5.92&  7.29& 5.02& 4.52& 5.31& 4.16& 3.38& 4.05& 3.92& 2.98& 3.26\\  
 15.0&  6.57&  7.88&  9.05& 5.96& 6.38&  7.02& 4.97& 4.47& 5.69& 3.96& 3.25& 3.61& 3.46& 2.83& 2.89\\  
 20.0&  6.74&  8.37& 10.81& 5.50& 6.57&  8.11& 4.95& 4.25& 4.99& 4.21& 3.84& 4.49& 3.60& 2.83& 2.97\\  
 25.0&  7.82&  9.53& 14.69& 6.21& 6.67&  7.09& 5.17& 4.28& 5.15& 4.02& 3.40& 3.83& 3.91& 2.95& 3.29\\  
 30.0&  6.50&  9.09& 10.46& 6.19& 6.82&  7.90& 4.68& 4.73& 4.91& 4.50& 3.67& 4.49& 3.90& 2.99& 3.13\\  
 35.0&  7.07&  8.87& 11.44& 6.33& 6.18&  8.74& 5.07& 4.72& 5.71& 4.27& 3.37& 4.39& 3.77& 3.02& 3.54\\  
 40.0&  7.53& 10.53& 13.78& 5.89& 6.91&  8.64& 4.55& 4.19& 5.09& 4.08& 3.34& 3.61& 3.89& 3.07& 3.08\\  
 50.0&  7.53& 10.16& 13.59& 6.08& 6.58&  7.12& 5.05& 5.29& 5.03& 4.32& 3.28& 4.03& 3.67& 3.00& 3.17\\
  \hline
\end{tabular}
\end{table}
\end{landscape}

\begin{landscape}
\begin{table}[!h]
\centering
\caption{Relative efficiencies of LASSO, aLASSO, and SCAD estimators when $n=100$, $r=0.9,$ $p=30$.} \label{tb:rel:eff:penalty:p30cor9}
\bigskip
\begin{tabular}{r|rrr|rrr|rrr|rrr|rrr}
\hline
 $\Delta^2$ & \multicolumn{3}{c|}{$k=1$} & \multicolumn{3}{c|}{$k=2$} & \multicolumn{3}{c|}{$k=3$}
 & \multicolumn{3}{c|}{$k=4$} & \multicolumn{3}{c}{$k=5$} \\
 & L & aL & SCAD & L & aL & SCAD & L & aL & SCAD & L & aL & SCAD & L & aL & SCAD \\
 \hline
   0.0& 29.79& 57.64& 11.91& 15.76& 18.35&  7.98& 14.37& 11.63& 7.51& 13.64& 9.22& 5.79& 10.57& 7.59& 4.06  \\  
   0.1& 24.18& 29.39& 14.37& 12.39& 13.05& 10.73& 15.68& 11.62& 7.97& 12.95& 9.02& 4.99&  8.38& 6.36& 3.43  \\  
   0.2& 14.87& 16.23& 12.63& 17.94& 17.19&  8.99& 15.24& 10.51& 7.40& 10.08& 7.93& 4.98&  9.32& 6.99& 3.76  \\  
   0.3& 26.00& 29.15& 15.69& 15.43& 18.23&  8.19& 13.82& 10.04& 6.24& 11.07& 8.59& 5.10&  9.53& 6.84& 3.67  \\  
   0.4& 28.22& 29.42& 14.06& 19.70& 15.41& 10.34& 12.88& 11.16& 7.28& 12.10& 8.52& 4.71& 10.25& 7.23& 3.75  \\  
   0.5& 22.06& 24.92& 13.13& 14.85& 14.17&  7.12& 14.33& 11.06& 6.36& 11.02& 7.24& 4.44&  8.92& 5.46& 3.48  \\  
   0.6& 25.28& 23.22& 14.83& 12.81& 12.69&  6.35& 13.34& 12.09& 6.05& 12.54& 6.91& 4.31&  8.87& 6.17& 3.50  \\  
   0.7& 18.98& 19.29&  9.90& 18.58& 17.69&  9.48& 13.58& 10.39& 5.90& 10.42& 6.60& 4.26&  9.73& 5.71& 3.67  \\  
   0.8& 20.48& 22.37&  9.36& 16.06& 15.32&  7.35& 14.41& 10.06& 6.13& 11.97& 7.09& 4.27&  8.65& 5.66& 3.72  \\  
   0.9& 17.53& 22.39& 10.94& 12.97& 12.74& 10.49& 11.49&  7.93& 4.99& 11.92& 7.84& 4.60&  8.40& 5.43& 3.21  \\  
   1.0& 18.85& 19.39& 16.86& 15.29& 12.60&  6.89& 13.28& 10.78& 6.53& 11.25& 7.66& 4.94&  9.07& 5.58& 3.07  \\  
   1.5& 18.81& 16.94&  8.83& 15.67& 12.90&  8.10& 12.70&  9.45& 5.94& 11.44& 7.38& 4.76&  8.34& 5.76& 3.27  \\  
   2.0& 18.73& 16.65& 11.90& 14.32& 10.13&  7.24& 10.07&  7.71& 3.99&  9.79& 6.45& 3.47&  8.17& 5.44& 3.37  \\  
   3.0& 13.71& 15.57& 11.84& 13.64&  9.60&  5.31& 10.15&  6.91& 4.23& 10.95& 6.10& 3.68&  8.11& 5.09& 2.98  \\  
   5.0& 13.90& 11.38&  6.54& 10.11&  8.23&  4.77&  9.73&  6.05& 3.90& 10.05& 6.31& 3.58&  8.93& 5.47& 3.19  \\  
  10.0& 10.88&  9.60&  6.02& 10.23&  7.41&  4.68&  8.78&  6.07& 3.68&  7.44& 4.97& 2.80&  7.38& 5.10& 2.93  \\  
  15.0&  9.81&  9.43&  5.00&  7.79&  6.68&  3.82& 10.43&  6.90& 4.40&  6.98& 4.54& 2.87&  7.26& 4.64& 2.94  \\  
  20.0& 12.36& 10.18&  6.50&  7.41&  6.79&  4.11&  9.02&  6.76& 4.17&  6.91& 5.05& 3.47&  7.70& 4.94& 2.84  \\  
  25.0&  9.93& 11.12&  6.48&  9.19& 10.29&  7.05&  8.42&  6.97& 4.28&  6.50& 4.97& 3.26&  6.46& 5.00& 3.08  \\  
  30.0&  9.53& 11.44&  7.09&  9.43&  7.97&  6.46&  7.81&  6.98& 4.82&  6.86& 5.73& 3.49&  6.65& 4.31& 3.20  \\  
  35.0& 10.19& 11.88&  7.00&  7.71&  8.33&  5.55&  7.53&  6.89& 4.70&  7.88& 5.91& 3.79&  7.21& 6.18& 4.36  \\  
  40.0&  9.54& 11.47&  7.21&  9.09& 10.96&  8.28&  8.11&  7.12& 4.99&  7.80& 7.04& 4.84&  5.83& 5.04& 3.43  \\  
  50.0& 10.31& 14.00&  9.59&  9.90& 13.83&  8.83&  8.70&  8.62& 6.10&  7.97& 6.89& 4.50&  6.98& 5.62& 3.85  \\
  \hline
\end{tabular}
\end{table}
\end{landscape}


\begin{landscape}
\begin{table}[!h]
\centering
\caption{Relative Efficiencies of LASSO, aLASSO, ElasticNet, and Stein-type Estimators when $r=0.2$, $n=100$, and $\Delta^2=0.$} \label{tb:penalty:JS:comparision}
\bigskip
\begin{tabular}{l|rrrrrrrrrrrr}
\hline
\backslashbox{Estimator}{$p$} & 10 & 15 & 20 & 25 & 30 & 40 & 50 & 60 & 70 & 80 & 90 & 95 \\
\hline
& \multicolumn{12}{c}{$k=0$} \\
\hline
S    &  4.79&  6.10&  7.58&  9.16& 11.89& 13.01& 12.73&  12.18&  10.49&   9.30&    4.24&     Inf \\
S+   &  7.22& 10.79& 14.03& 16.42& 23.42& 24.60& 23.06&  27.02&  21.78&  19.97&    9.62&     Inf \\
L    & 10.51& 20.26& 18.95& 26.97& 29.98& 54.79& 61.65&  72.87& 151.25& 233.12&  617.51& 1282.98 \\
aL   & 12.30& 23.33& 37.17& 33.39& 54.50& 98.54& 90.24& 156.91& 249.65& 405.12& 1175.86& 2395.63 \\
EN25&  11.94& 17.06& 26.58& 30.51& 38.42& 56.08& 75.12& 103.92& 197.13& 312.53& 665.46& 2260.51 \\
EN50&  10.90& 16.26& 21.93& 28.55& 34.78& 51.86& 71.75&  82.26& 146.58& 308.54& 616.73& 1783.79 \\
EN75&  10.00& 14.66& 21.65& 25.86& 33.07& 41.84& 61.88&  75.16& 143.15& 279.61& 608.76& 1844.96 \\
\hline
& \multicolumn{12}{c}{$k=1$} \\
\hline
S   & 3.89& 5.88&  6.89&  8.51&  9.96& 11.29& 12.42& 10.64& 10.35&   8.45&   4.42&    Inf \\
S+  & 5.37& 8.50& 10.51& 13.66& 15.59& 19.17& 22.78& 22.27& 20.71&  17.70&   9.26&    Inf \\
L   & 3.89& 6.14&  8.44&  9.65& 12.72& 18.82& 24.23& 44.00& 58.47&  90.62& 301.28& 648.22 \\
aL  & 3.97& 5.83&  8.86& 11.36& 12.27& 19.59& 28.18& 59.66& 85.39& 149.74& 325.10& 869.17 \\
EN25&  4.22& 6.48& 9.09& 30.51& 14.22& 23.43& 31.74& 47.75& 73.15& 133.89& 334.52&  997.93 \\
EN50&  4.10& 6.23& 8.64& 28.55& 13.85& 22.74& 30.62& 43.62& 74.54& 129.17& 283.49& 1018.99 \\
EN75&  3.91& 6.28& 8.56& 25.86& 13.73& 22.53& 31.02& 40.82& 68.85& 122.62& 288.95&  925.23 \\
\hline
& \multicolumn{12}{c}{$k=3$} \\
\hline
S  & 2.99& 4.29& 5.81& 6.96& 8.06&  9.57& 11.11& 10.64& 10.02&  7.36&   4.13&    Inf \\
S+ & 3.35& 4.83& 6.68& 8.18& 9.83& 12.94& 14.82& 15.95& 15.60& 13.28&   9.41&    Inf \\
L  & 2.06& 3.00& 4.11& 5.18& 6.71&  9.16& 13.41& 22.13& 31.50& 57.30& 152.62& 378.14 \\
aL & 1.56& 2.40& 3.02& 4.00& 5.29&  7.85& 11.46& 19.03& 29.85& 48.55& 127.24& 382.55 \\
EN25&  2.39& 3.35& 4.46& 5.76& 7.32& 11.12& 16.29& 23.88& 36.54& 63.86& 162.73& 441.96 \\
EN50&  2.20& 3.16& 4.11& 5.61& 6.96& 10.56& 15.60& 22.91& 34.76& 60.21& 150.59& 424.36 \\
EN75&  2.14& 3.02& 4.07& 5.46& 6.76& 10.12& 15.19& 21.75& 33.87& 58.47& 146.62& 403.01 \\
\hline
& \multicolumn{12}{c}{$k=5$} \\
\hline
S  & 2.35& 3.38& 4.34& 5.33& 6.26& 7.95&  9.59& 10.17& 10.22&  8.33&  4.50&    Inf \\
S+ & 2.39& 3.48& 4.57& 5.64& 6.63& 8.95& 11.06& 12.71& 13.59& 12.80&  7.72&    Inf \\
L  & 1.46& 2.14& 2.80& 3.63& 4.45& 6.90&  9.82& 14.04& 22.64& 35.90& 90.95& 286.14 \\
aL & 1.09& 1.56& 1.96& 2.50& 3.08& 4.70&  6.83& 10.33& 16.40& 26.74& 71.77& 219.63 \\
EN25  &1.79& 2.50& 4.46& 4.11& 5.20& 11.12& 10.68& 16.41& 36.54& 44.92& 108.73& 441.96 \\
EN50  &1.63& 2.35& 4.11& 3.81& 4.82& 10.56&  9.93& 15.23& 34.76& 40.24& 100.50& 424.36 \\
EN75  &1.55& 2.24& 4.07& 3.69& 4.68& 10.12&  9.49& 14.83& 33.87& 40.20&  98.05& 403.01 \\
\hline
\end{tabular}
\end{table}
\end{landscape}

\bibliographystyle{apa}
\bibliography{ShrinkageReferences}
\renewcommand{\leftmark}{References} 

 \end{document}